\documentclass[twocolumn]{autart}    

\usepackage{xcolor}
\usepackage{graphicx} 
\usepackage{amsmath} 
\usepackage{amssymb}  
\usepackage{dsfont}
\usepackage{algpseudocode}
\usepackage{algorithm}
\usepackage[scr]{rsfso}
\usepackage[percent]{overpic}
\usepackage{cite} 

\DeclareMathOperator*{\argmin}{arg\,min}
\usepackage{caption}
\usepackage{subcaption}
\usepackage{color}
\usepackage{enumitem}
\usepackage{mathtools}

\newtheorem{definition}{Definition}
\newtheorem{assumption}{Assumption}

\newtheorem{lemma}{Lemma}

\newtheorem{remark}{Remark}
\newcommand{\R}{\mathbb{R}}
\newcommand{\xminus}{x_{-i}}
\newcommand{\viol}{\mathcal{V}}
\newcommand{\violset}{\mathbb{V}}
\newcommand{\algo}{\mathscr{A}}
\newcommand{\trasp}{\ensuremath{\raisebox{1pt}{$^{\intercal}$}}}
\newcommand{\traspm}{\ensuremath{^{\intercal}}}
\newcommand{\msample}{\boldsymbol{\delta}_K}
\newcommand{\ball}[1][x^*]{\ensuremath{\mathbb{B}}({#1}, \rho)}
\newcommand{\vtight}{\boldsymbol{\rho}}
\newcommand{\tight}{c\rho}
\newcommand{\rvline}{\hspace*{-\arraycolsep}\vline\hspace*{-\arraycolsep}}
\newcommand{\bigzero}{\mbox{\normalfont\large\bfseries 0}}
\newcommand{\bfzero}{\mathbf{0}}
\newcommand{\vQ}{\ensuremath{\vec{Q}}}
\newcommand{\mukappa}{\mu^{(\kappa)}}
\newcommand{\erre}[1][\mu]{\ensuremath{\mathcal{L}_{#1}}}
\newcommand{\setX}{\ensuremath{X}}
\newcommand{\setN}{\ensuremath{\mathcal{N}}}
\newcommand{\setM}{\ensuremath{\mathcal{M}}}
\newcommand{\setW}{\ensuremath{\mathcal{U}}}

\newcommand{\Dsqr}[1][-]{\ensuremath{D_s^{{#1}1/2}}}
\newcommand{\TD}{\ensuremath{\tilde{T}_D}}
\newcommand{\uppb}{\ensuremath{\overline{b}_{\delta_k}}}
\newcommand{\lowb}{\ensuremath{\underline{b}_{\delta_k}}}
\newcommand{\Dmu}{\ensuremath{D_{\mu}^{-1}}}
\usepackage{pifont}
\newcommand{\cmark}{\ding{51}}%
\newcommand{\xmark}{\ding{55}}%

\newcommand{\KM}[1]{{\color{black} #1}}
\newcommand{\GP}[1]{{\color{black} #1}}
\newcommand{\FF}[1]{{\color{black} #1}}

\begin{document}

\begin{frontmatter}
\title{A priori data-driven robustness guarantees on strategic deviations from generalised Nash equilibria}
\author[TUD]{George Pantazis}\ead{G.Pantazis@tudelft.nl},    
\author[US]{Filiberto Fele}\ead{ffele@us.es},               
\author[UOXF]{Kostas Margellos}\ead{kostas.margellos@eng.ox.ac.uk}  

\address[TUD]{Delft Center for Systems and Control, TU Delft, 2628 CD Delft, Netherlands}  
\address[US]{Department of Systems Engineering and Automation, University of Sevilla, Spain}  
\address[UOXF]{Department of Engineering Science, University of Oxford, OX1 3PJ, UK}  

\begin{keyword}                           
Generalised equilibrium problem;  Randomized methods; Stochastic game theory; Multi-agent systems. 
\end{keyword}                             

\begin{abstract}                         
In this paper we focus on noncooperative games with uncertain constraints coupling the agents' decisions. We consider a setting where bounded deviations of agents' decisions from the equilibrium are possible, and uncertain constraints are inferred from data. Building upon recent advances in the so called scenario approach, we propose a randomised algorithm that returns a nominal equilibrium such that a \emph{pre-specified} bound on the probability of violation for yet unseen constraints is satisfied for an entire region of admissible deviations surrounding it---thus supporting neighbourhoods of equilibria with probabilistic feasibility certificates. 
For the case in which the game admits a potential function, whose minimum coincides with the social welfare optimum of the population, the proposed algorithmic scheme opens the road to achieve a trade-off between the guaranteed feasibility levels of the region surrounding the nominal equilibrium, and its system-level efficiency. Detailed numerical simulations corroborate our theoretical results. 
\end{abstract}

\end{frontmatter}

\section{Introduction}
The study of noncooperative games plays a significant role in a panoply of applications ranging from smart-grids \cite{Saad1} to communication \cite{ScutaryMonotone} and social networks \cite{Acemoglu}. In these setups, agents can be modelled as self-interested entities that interact with each other and make decisions based on possibly misaligned individual criteria, while being subject to constraints (local or global) that restrict the domain of their choices.
Even though a variety of systems can be analysed by means of deterministic game-theoretic tools \cite{ScutaryMonotone,GrammaticoEtAl2016,Dario2019}, in many applications the decision making procedure is affected by uncertainty. 
A number of results in the literature have explicitly addressed uncertainty in a noncooperative setting. Specifically, \cite{Bopardikar_2013}  follows a randomized approach for the special case of stochastic zero-sum games. Most results rely on specific assumptions on the probability distribution \cite{Kouvaritakis,Singh} and/or the geometry of the uncertainty set \cite{Aghassi2006,FukuSOCCP,Nishimura_2009}.

To circumvent these limitations, recent developments adopt a data-driven perspective, focusing on the connection of game theory with the so called scenario approach  \cite{CalafioreCampi2006}. This is based on the idea that an optimisation problem with constraints parametrised by an uncertain parameter---with fixed but possibly unknown support set and probability distribution---can be approximated by drawing samples of that parameter and solving the problem subject to the constraints generated by those samples only; this approximation is known as the \emph{scenario program}. Standard results in the scenario approach \cite{CampiGaratti2008,CampiEtAl2009,Calafiore_2010} provide certificates on the probability that a new yet unseen constraint will violate the randomised solution obtained by the scenario program. 

While the aforementioned results apply to uncertain \emph{convex} optimisation problems, the works \cite{CampiGaratti2018_waitjudge} and \cite{CampiGarattiRamponi2018} paved the way towards the provision of data-driven robustness guarantees to solutions of more general nonconvex problems. 
In \cite{Feleconf2019,Fele2021,Dario_Scenario}, these theoretical advancements were leveraged for the first time in a game-theoretic context, for the formulation of distribution-free probabilistic feasibility guarantees for randomised Nash equilibria. These works provide guarantees for one specific equilibrium point (often assumed to be unique); this was extended in \cite{fele-a,Pantazis2020}, by providing \emph{a posteriori} feasibility guarantees for the entire domain. Besides the game-theoretic context, alternative methodologies for set-oriented probabilistic feasibility guarantees have been proposed in the seminal works \cite{DeFariasDP_24,Alamo2009_11}, which a priori characterise probabilistic feasibility regions constructed out of sampled constraints using statistical learning theoretic results. More recently, the so called probabilistic scaling \cite{Alamo_2019_14,Mammarella2022_15} has been proposed to obtain a posteriori guarantees on the probability that a polytope generated out of samples is a subset of some chance-constrained feasibility region. Following an approach similar to \cite{fele-a}, the works \cite{Fabiani2020a,Fabiani2020b} deliver tighter guarantees by focusing on variational-inequality (VI) solution sets.

The results above follow a standard approach in the game-theoretic literature, where a strict behavioural assumption---the so called \emph{rationality}---is imposed on the players' decision making. Namely, the players are viewed as rational agents wishing to maximize their profit (expressed by some given cost function). However, studies have shown that this is unrealistic in practice \cite{Vernon2008,Camerer,Prospect1,Prospect4} and that agents usually exhibit a \emph{boundedly rational} behaviour \cite{Rubinstein1997}, i.e., their decisions can deviate from rationality due to individual biases, behavioural inertia, restricted computational power/time, etc. 
The consequences of this become relevant in engineering applications, as the human role in technical systems evolves beyond mere users and consumers to active agents, operators, decision-makers and enablers of efficient, resilient and sustainable infrastructures \cite{LAMNABHILAGARRIGUE2017}. 

To bridge this gap between real-world applications and the cognate literature, here we study games with uncertain constraints, where deviations from a \emph{nominal} equilibrium are explicitly considered. We follow a randomised approach to approximate the coupling constraints by means of data. In this more general setting, where deviations are considered, providing guarantees for a single solution is devoid of any meaning: indeed, repetition of the game might lead to a different solution in a neighbourhood around the nominal equilibrium, irrespective of the employed dataset. Technically speaking, this renders the identification  of the data samples that support the solution (cf.~sample compression \cite{MargellosEtAl2015}) a challenging task. 
Focusing on the class of generalised Nash/Wardrop equilibrium seeking problems (GNEs/GWEs)  \cite{FacchineiKanzow2009}, we contribute to the provision of data-driven robustness guarantees for the collection of possible deviations from the equilibrium as follows:
\begin{enumerate}
\item  \KM{Adopting a scenario-theoretic paradigm, we establish a methodology for the provision of \emph{a posteriori} probabilistic feasibility guarantees for a region around the randomised equilibrium of the game under study. This result (Theorem 1), complements \cite{fele-a}, \cite{Pantazis2020}, \cite{Fabiani2020a}, \cite{Fabiani2020b}, that \FF{instead} focus on the entire feasibility region. Focusing on a \FF{circumscribed region around a GNE/GWE allows offering tighter} probabilistic bounds, while the results of \cite{fele-a}, \cite{Pantazis2020}, \cite{Fabiani2020a}, \cite{Fabiani2020b}, can be obtained as a limiting case of Theorem 1.}
\item \KM{We design a data-driven algorithm that returns a GNE/GWE and show that all points
in a predefined admissible region surrounding it enjoy \emph{a priori} probabilistic feasibility guarantees. This result (Theorems 2 and 3), unlike Theorem 1, offers an \emph{a priori} statement \FF{valid for a region that is tunable} by the user, modelling possible deviations from a nominal equilibrium that a designer wishes to take into account when incentivising a certain operation profile. \\
A distinctive feature of this result is that it provides \emph{a priori} guarantees for a set rather than \FF{single points} \cite{Dario_Scenario}, \cite{Feleconf2019}, \cite{Fele2021}. \FF{These} guarantees depend on a prespecified quantity, which in turn can affect the location of the nominal equilibrium and the size of the region for which these probabilistic guarantees hold. As such this region is tunable, unlike \cite{GrammaticoNonConvex} where \emph{a priori} guarantees for a set of solutions are provided, but this set is not controlled by the user and could be arbitrarily narrow. Moreover, the results of \cite{GrammaticoNonConvex}  do not focus on games and follow a fundamentally different approach.}
\\ 
Furthermore, when the game under study admits a potential function---whose minimum coincides with some social welfare optimum---our methodology provides a new perspective for trading off the probabilistic feasibility of the region surrounding the nominal equilibrium and its system-level efficiency.
\item \FF{We propose an equilibrium seeking algorithm as the machinery to obtain a region surrounding a GNE/GWE over which the aforementioned feasibility guarantees hold. The algorithm relies on a primal-dual scheme and is} \KM{inspired by seminal developments in \cite{Pang1}. However, the mapping that characterizes the algorithm updates differs from those typically encountered in the literature (e.g., see \cite[Ch.~12]{Pang1}). This requires showing that the ad-hoc mapping enjoys certain continuity and co-coercivity properties, thus extending the proof-line of \cite{Pang1} (see Lemmas 2 \& 3, and proof of Theorem 2), a task which is interesting per se.
}
\end{enumerate}

\setlength{\arrayrulewidth}{0.2mm}
\setlength{\tabcolsep}{10pt}
\renewcommand{\arraystretch}{1.3}
\begin{table*}[h] 
	\centering
\begin{tabular}{ |p{4cm}|p{4.5cm}|p{2.8cm}|p{0.7cm}| p{2cm}|  }
	\multicolumn{5}{c}{Table 1 - Classification of results and comparison with cognate literature} \\
	\hline
	Problem class & Solution sets  &Type of feasibility guarantees & Tuning & References  \\
	\hline
	\hline
			Affine feasibility problems   & Entire feasible set    &a posteriori &   ~~~\xmark   &  \hspace{-0.02cm}\cite[Thm. 6]{fele-a} \\
		Convex feasibility problems &   Entire feasible set  & a posteriori   &  ~~~\xmark  & \hspace{-0.02cm}\cite[Thm. 2]{Pantazis2020} \\ 
				Uncertain GNEs & GNE solution set & a posteriori &   ~~~\xmark   & \hspace{-0.02cm}\cite[Thm. 1]{Fabiani2020a}, \cite[Thm. 1]{Fabiani2020b}\\
\textbf{Uncertain GNEs/GWEs} &   \textbf{Subset of feasible deviations around GNE/GWE} & \textbf{a posteriori}   &  ~~~\xmark  & \textbf{Theorem 1} \\ 
	Uncertain VIs  & Unique solution    &a priori/a posteriori & ~~~\xmark  &\hspace{-0.02cm}\cite[Cor. 1]{Dario_Scenario}, \cite[Thm. 5]{Feleconf2019},\cite[Thm. 8]{Fele2021}  \\
		Convex feasibility problems &   (Arbitrary) inner approx. of feasible set  &  a priori    &  ~~~\xmark  &  \hspace{-0.02cm}\cite[Thm.2]{GrammaticoNonConvex} \\
		\textbf{Uncertain GNEs/GWEs} &   \textbf{Tunable subset of feasible deviations around GNE/GWE} & \textbf{a priori}   &  ~~~\cmark  & \textbf{Theorem 3} \\ 
	\hline
\end{tabular}
\end{table*}

Our contributions compared to the cognate literature are summarized in Table 1. The rest of the paper is organized as follows. In Section \ref{sec_prelim} we provide fundamentals of game theory and the scenario approach. In Section \ref{sec_aposteriori} we show how the feasibility guarantees for a region around the game solution can be a posteriori quantified. In Section~\ref{sec_apriori} we propose a data driven algorithm and prove its convergence to an equilibrium such that the considered neighbourhood of strategic deviations can satisfy prespecified probabilistic feasibility requirements. An illustrative example in Section \ref{sec:num_example} corroborates our theoretical analysis. Section \ref{sec:conclusion} concludes the paper and presents future research directions. To streamline the presentation of our results, some proofs are deferred to the Appendix.

\section{Preliminaries} \label{sec_prelim}
\emph{Notation}: All vectors are column unless otherwise indicated. $\mathbb{R}^n_{+}$ is the nonnegative orthant in $\mathbb{R}^n$. For an $n\times n$ matrix $A$, we write $A \succ 0$ ($A \succeq 0$) when $x\trasp A x > 0$ ($x\trasp A x \geq 0$), for any $x\in\mathbb{R}^n$.
We denote by $\mathbf{0}_{q \times r}$ the $q \times r$ null matrix, by $I_r$ the $r \times r$ identity matrix, and by $\mathbf{1}_r$ the vector of $r$ ones; dimensions can be omitted when clear from the context. $e_{q}$ is the unit vector whose $q$-th element is 1 and all other elements are 0, $\|\cdot\|_p$ the $p$-norm operator, and  $(\cdot)_r$ denotes the $r$-th component of its vector argument. $\mathbb{B}_p(x,\rho)=\{y \in \mathbb{R}^{d}: \|y-x\|_p < \rho\}$ is the open $p$-normed ball centred at $x$ with radius $\rho$; when $p$ is omitted, any choice of norm is valid. For a set $S$, $|S|$ denotes its cardinality, while $2^S$ denotes its power set, i.e., the collection of all subsets of $S$. Finally, given $D\succ 0$, $\text{proj}_{K,D}[x]:= \argmin_{y \in K}(y-x)\trasp D(y-x)$ is the skewed projection of $x$ onto the set $K$.

\subsection{Games with uncertain constraints}
\label{sec:aggr_games}
We consider a population of agents  with index set $\mathcal{N}=\{1, \dots, N\}$. The decision vector $x_i$ of each agent $i \in \mathcal{N}$ takes value in the set $X_i \subseteq \R^{n}$, while $x=(x_i)_{i=1}^N \in \setX =\prod_{i=1}^N X_i \subseteq  \R^{nN}$ is the global decision vector that is formed by concatenating the decisions of the entire population. The vector $\xminus \in \R^{n(N-1)}$ comprises all agents' decisions except for those of agent $i$. In our setup,  the cost incurred by agent $i \in \mathcal{N}$ is expressed by a real-valued function $J_i(x_i, \xminus)$ that depends on local decisions as well as on the decisions from other agents $j\in\setN\setminus\{i\}$. In the following, with a slight abuse of notation, we can exchange $x$ for $(x_i,\xminus)$ to single out agent $i$'s decision from the global decision vector.
Furthermore, we consider \emph{uncertain} constraints coupling the agents' decisions. These can be expressed in the form\footnote{This formulation can describe deterministic and/or local constraints as special cases.}
\begin{equation}
	C_\delta=\{x \in X: g(x,\delta) \leq 0 \},\; \delta\in\Delta, \label{constraint_delta}
\end{equation}
where $g: \R^{nN} \times \Delta \rightarrow \R$ depends on some uncertain parameter $\delta $ taking values in a support set $\Delta$ according to a probability measure $\mathbb{P}$.

Feasible collective decisions under this setup can be found by letting every agent $i \in \mathcal{N}$ solve the following optimization program, where the decisions $\xminus$ of all other agents are given, 
\begin{equation}\label{eq:game}
	\left.
	\begin{aligned}
		G : \;&\min_{x_i \in X_i} J_i(x_i, \xminus) \\ 
		&\text{ subject to } x_{i} \in \bigcap_{\delta \in \Delta} C^i_{\delta}(\xminus)  
	\end{aligned}
	\right\}\; \forall i\in\mathcal{N};
\end{equation}
here, $C^i_{\delta}(\xminus)= \{ x_i \in X_i: g(x_i, \xminus, \delta ) \leq 0 \}$ is the projection of the coupling constraint on $X_i$ for fixed $\xminus$ and uncertain realization $\delta \in \Delta$.  The collection of coupled optimization programs in \eqref{eq:game} for all $i \in \mathcal{N}$ constitutes an \emph{uncertain noncooperative game}; we denote it as $G$. 

Note that \eqref{eq:game} follows a worst-case paradigm, taking into account all possible coupling constraints that can be realised by variations of the uncertain parameter $\delta \in \Delta$. 
This can render the solutions of $G$ rather conservative. Furthermore, it is in general not possible to compute a solution for $G$ without an accurate knowledge of, and/or  additional assumptions on, the support set $\Delta$ and the probability distribution $\mathbb{P}$. To circumvent these limitations, we follow a data-driven paradigm and approximate  $G$ by means of a finite number of samples drawn from $\Delta$, namely the $K$-multisample $\msample=(\delta_1, \dots, \delta_K)\in\Delta^K$. In the sequel, we hold on to the standing assumption that these samples are independent and identically distributed (i.i.d.). Apart from this, no other knowledge on the support set $\Delta$ and the probability distribution $\mathbb{P}$ of the uncertain parameter is required.
Then, for a given multi-sample $\msample$, \eqref{eq:game} can be rewritten as
\begin{equation}\label{eq:scenario_game}
	\left.
	\begin{aligned}
	G_K :\; &\min_{x_i \in X_i} J_i(x_i, \xminus) \\  
	& \text{ subject to } \ x_{i} \in \bigcap_{k=1}^K C^i_{\delta_k}(\xminus)
	\end{aligned}
	\right\}\; \forall i\in\mathcal{N}.
\end{equation}
Instead of considering all possible uncertainty realizations $\delta \in \Delta$ as in \eqref{eq:game}, we let the data encoded in $\msample$ lead agents to their decision by solving \eqref{eq:scenario_game}. We refer to the collection of coupled optimization programs in \eqref{eq:scenario_game} as the \emph{scenario game} $G_K$ (the subscript $K$ implies dependence on the drawn multi-sample $\msample$). Under standard assumptions, a solution to the scenario game $G_K$ exists and the problem is, in contrast to $G$, tractable using state-of-the-art equilibrium seeking algorithms.  

\subsection{Variational inequalities and game equilibria}
\label{sec_gamesolVI}
Notably---under certain assumptions detailed next---solutions to the game $G_K$ can be retrieved as solutions to a variational inequality (VI),  for specific choices of the mapping $F:\, X\to\mathbb{R}^{nN}$ \cite[Thm~3.9]{FacchineiKanzow2009}:
\begin{equation}\label{eq:VI_K}
	\begin{split}
		\text{VI}_K:& \ \text{Find} \ x^* \in \Pi_K  \ \text{such that} \\ 
		& \ (x- x^*)\trasp F(x^*) \geq 0 \ \text{for any} \ x \in \Pi_K,
	\end{split}
\end{equation}
where $\Pi_K:=X \cap \bigcap_{k=1}^{K} C_{\delta_k}$ denotes the problem domain.
A classic game solution concept, which encounters wide application in the literature, is the Nash equilibrium (NE) \cite{nash50}. At a NE, no agent can decrease their cost by unilaterally changing their decision. Formally, this can be stated as follows.
\begin{definition} \label{def_GNE}
	A point $x^*=(x_i^*,\xminus^*)\in \Pi_K$ is called a generalised Nash equilibrium (GNE) of $G_K$ if, for all $i \in \setN$,
	\vspace{-1em}
	\begin{equation*}
		J_i \left(x_{i}^*,\xminus^* \right) \leq J_i (y_i, \xminus^*),\; \forall y_i \in X_i \cap \bigcap_{k=1}^K C^i_{\delta_k}(\xminus^*).
	\end{equation*}
\end{definition}

For our analysis, we rely on the following conditions:
\begin{assumption}  \label{assum_cvx_differ} 
	For all $i \in \mathcal{N}$, $J_i(x_i,\xminus)$ is convex and continuously differentiable in $x_i$ for any fixed $\xminus$.
\end{assumption}
\begin{assumption} \label{assum_affine}
	\begin{enumerate}
		\item For any multi-sample $\msample \in \Delta^K $, the domain $\Pi_K$ is non-empty.
		\item The set $X= \prod_{i=1}^N X_i$ is compact, polytopic and convex.
		\item For any $\delta \in \Delta$,  $g$ is an affine function of the form $g(x, \delta)=a(\delta)\trasp x -b(\delta)$, where $a: \Delta \rightarrow \R^{nN}$ and $b: \Delta \rightarrow \R$. 
	\end{enumerate}
\end{assumption}
 Note that convexity of the cost function with respect to the agent's local decision is crucial for the design of tailored algorithms with theoretical convergence guarantees for Nash equilibrium seeking.
Under these assumptions, we can determine a GNE as in Definition~\ref{def_GNE} by solving \eqref{eq:VI_K} with 
\begin{equation}\label{eq:F}
	F(x)=F_{\mathrm{NE}}(x):=\begin{bmatrix}
		\nabla_{x_1}J_1(x_1,x_{-1})\\
		\vdots\\
		\nabla_{x_N}J_N(x_N,x_{-N})
	\end{bmatrix}.
\end{equation}
A class of problems of common interest can be modelled by the so called \emph{aggregative} games \cite{KUKUSHKIN2004,Jensen2010,Acemoglu2013}, where the cost incurred by agents depends on some aggregate measure---typically the average---of the decision of the entire population. Such a cost can be expressed in \eqref{eq:scenario_game} by the function $J_i(x_i, \sigma(x))$, where the aggregate $\sigma: \R^{nN} \rightarrow \R^n$ is defined as the mapping $x \mapsto \frac{1}{N}\sum_{i=1}^N x_i$. A solution  frequently linked to this class of games is the Wardrop equilibrium (WE), a concept akin to the NE but specifically defined in the context of transportation networks \cite{BeckmannEtAlBOOK}. The variational WEs of $G_K$ can be expressed by using $F(x)=F_{\mathrm{WE}}(x):=[\nabla_{x_i}J_i(x_i,z)_{|z=\sigma(x)}]_{i \in \mathcal{N}}$; notice that in this case the second argument of $J_i$ is fixed and set to $\sigma(x)$, consistently with the notion of WE where agents neglect the impact of their decision on others.

We restrict the considered class of variational mappings as follows:
\begin{assumption} \label{assum_monotonicity_continuity}
	The mapping $F$ is
	\begin{enumerate}
		\item  $\alpha$-strongly monotone, i.e., $(x-y)\trasp(F(x)-F(y)) \geq \alpha \|x-y\|^2$  for any  $x, y \in X$, \label{strong_monotonicity}
		\item  $L_F$-Lipschitz continuous., i.e., 
		$\|F(x)-F(y)\|  \leq L_F\|x-y\|$ for any $x, y \in X$. \label{Lipschitz}
	\end{enumerate}
\end{assumption}
Assumptions \ref{assum_cvx_differ} and \ref{assum_monotonicity_continuity} are standard in the game-theoretic literature \cite{Pang1,ScutaryMonotone}. Assumption \ref{assum_affine} is relatively mild; the affine form of the constraints is exploited in the proposed algorithm (see Section 3) for the convergence to an equilibrium bearing the desired robustness properties. 

We point out that in general only a subset of solutions to $G_K$ can be retrieved through \eqref{eq:VI_K}: these are referred to as \emph{variational equilibria}, and enjoy favourable properties over nonvariational ones, as with the former the coupling constraints' burden is equally split among agents \cite{HARKER1991,KULKARNI2012}. 
The following lemma, adapted from \cite[Thm.~2.3.3]{Pang1}, formalises the connection between the solutions to VI$_K$ and the GNEs (or GWEs) of $G_K$.
\begin{lemma}  \label{existence}
	Under Assumptions \ref{assum_cvx_differ}, \ref{assum_affine} and \ref{assum_monotonicity_continuity}, 
	VI$_K$ has a unique solution that is also an equilibrium of $G_K$.
\end{lemma}

For the considered  class of VIs, several algorithms from the literature can be employed to obtain a variational equilibrium of $G_K$; see, e.g., \cite{FacchineiKanzow2009,Dario2019}. We remark that, even if not explicitly shown for ease of notation, any solution $x^*$ to $G_K$ is itself a function of the drawn multisample $\msample\in\Delta^K$.
Probabilistic feasibility guarantees for the unique solution of VI$_K$ can then be provided both in an \emph{a priori} and \emph{a posteriori} fashion by resorting to the results in \cite{Feleconf2019,Fele2021,Dario_Scenario}. However, these results are tailored to the provision of probabilistic feasibility guarantees for a single point (namely the solution of a VI): any strategic deviation from the equilibrium is not covered by such guarantees. We cover this issue in Section~\ref{sec:probfeas_sets}. First, we provide some background on the scenario approach.

\subsection{Basic concepts in the scenario approach}\label{sec:background_scenario_approach}
A fundamental notion in the scenario approach is the \emph{probability of violation} of an uncertain constraint.
\begin{definition} \label{violation} 
	\begin{enumerate}
		\item 	The probability of violation $\viol: \mathbb{R}^{nN} \rightarrow [0,1]$ of a point $x \in \Pi_K$ is defined as the probability that a new yet unseen sample $\delta \in \Delta$ will give rise to a constraint $C_\delta$ (as defined in \eqref{constraint_delta}) such that $x \notin C_\delta$, i.e.,  
	$\viol(x):=\mathbb{P}\{\delta \in \Delta: x \notin C_\delta\}.$
		\item The probability of violation $\violset: 2^{\mathbb{R}^{nN}} \rightarrow [0,1]$ of a set $S\subseteq \Pi_K$ is defined as the worst-case $\viol$ among all the points in $S$, i.e., $\violset(S)=\sup\limits_{x \in S}\mathbb{P}\{\delta \in \Delta: x \notin C_\delta\}$.
	\end{enumerate} 
\end{definition}
A data-driven decision-making process can be formally characterized by a mapping---the \emph{algorithm}---that takes as input the data encoded by the samples and returns a set of decisions.
\begin{definition} \label{Algorithm_definition}
An algorithm is a function $\algo: \Delta^{l} \rightarrow 2^{\mathbb{R}^{nN}} \times \mathbb{R}^{nN}$ that takes as input an $l$-multisample and returns the \FF{pair $(x^\ast, S_l^*(x^\ast))$, namely,  a point $x^\ast$ and a solution set $S_l^\ast$ parametrized by  $x^\ast$.}
\end{definition}
\FF{In the setting we consider, we have $x^* \in S_l^*(x^*)$; however, this ought not to be the case in general.}
In the following, we interpret the above definition as context-dependent, in that the size $l$ of the input multisample is admitted to vary---all else remaining fixed for a given algorithm $\algo$. 

A key notion, strongly linked to that of algorithm, is the \emph{minimal compression set} \cite{MargellosEtAl2015}. This concept springs from the observation that typically only a subset of the sampled data is relevant to a decision or set of decisions, and all other samples are redundant.

\begin{definition}[Compression set] \label{def_compression}
	Consider an algorithm $\algo$ as in Definition~\ref{Algorithm_definition}. A subset of samples $I \subseteq \msample$  is called a \emph{compression} for $\algo(\msample)$ if  $\algo(I) = \algo(\msample)$.\footnote{With some abuse of notation, in the remainder the symbol $\msample$ is interpreted as either the i.i.d.~sample vector $\msample\in\Delta^K$, or the set comprising its components, i.e., $\msample=\{\delta_1,\ldots,\delta_K\}\subseteq \Delta$, depending on the context.} 
	As multiple subset of samples can exist that fulfil this property, the ones with the minimal cardinality are called \emph{minimal} compression sets.
\end{definition} 
If we feed the algorithm with the set of samples corresponding to a compression, then the same decision will be returned as when we feed the algorithm with the entire multi-sample. As established in \cite{MargellosEtAl2015}, the compression set is related to the notion of support samples \cite{CampiGaratti2008} and that of essential constraints \cite{Calafiore_2010}. Under certain non-degeneracy assumptions these concepts coincide.



\section{Probabilistic feasibility of sets around equilibria}
\label{sec:probfeas_sets}
\subsection{A first a posteriori result}\label{sec_aposteriori}
Returning to the scenario game $G_K$ in \eqref{eq:scenario_game}, we now consider a more general setup where agents are allowed to deviate from  $x^*$ following, e.g., unmodelled changes in their cost functions; while we suppose that these deviations are feasible with respect to the local constraints, we want to study the feasibility as regards the coupling constraints obtained through sampling. Specifically, the region in which agents' strategies can deviate from the nominal equilibrium is assumed to lie within a predefined open ball $\mathbb{B}(x^*, \rho)$, where $\rho>0$ is a fixed radius that denotes the maximum possible distance of agents' deviations from $x^*$; the latter is assumed to be unique as per Lemma~\ref{existence}. As such, the region of interest is $S_K^*=\Pi_K \cap \mathbb{B}(x^*, \rho)$.  

\begin{figure}
	\centering
	\begin{overpic}[scale=0.5,trim=6.5cm 4.55cm 2.5cm 4.7cm,clip]{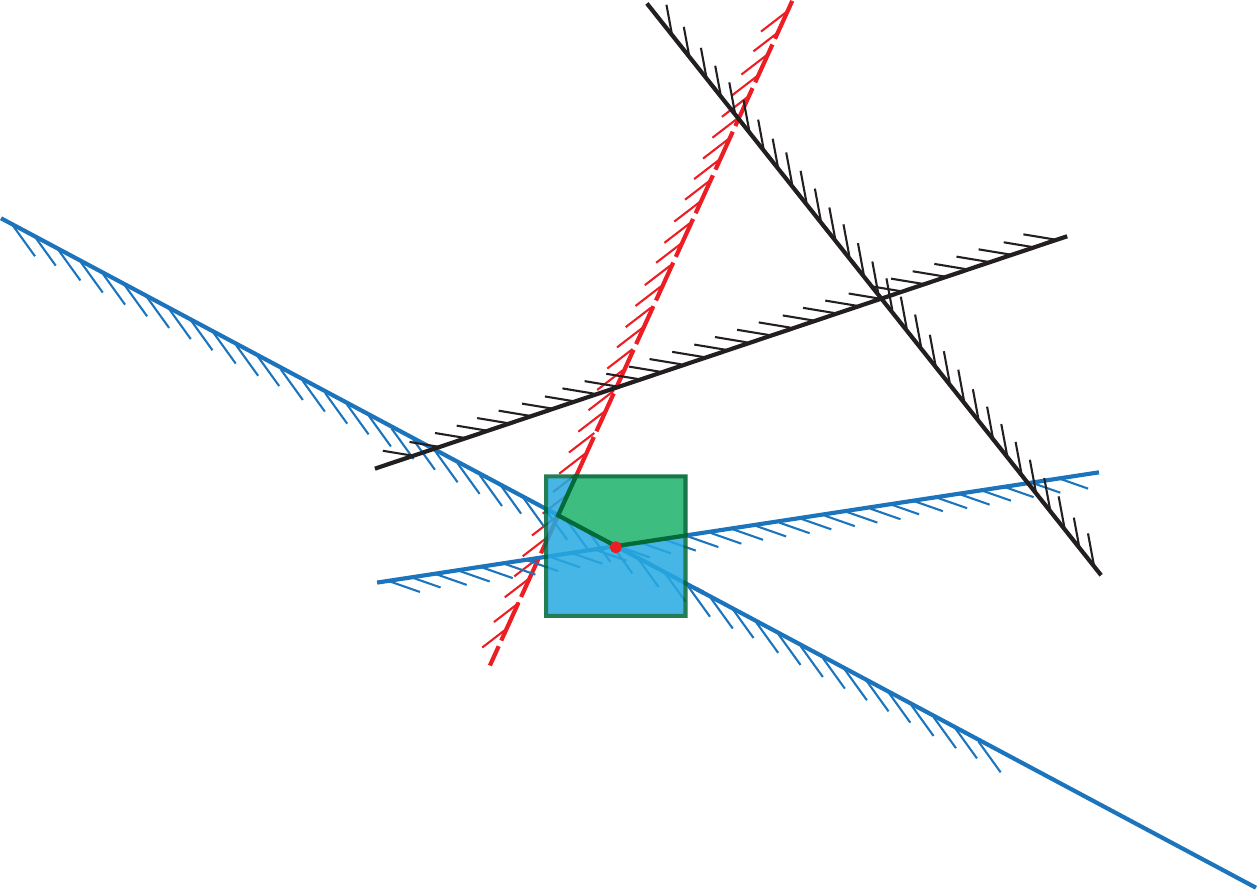}
		\put (30,4) {$x^\ast$}
		\put (47,25) {$\Pi_K$}
		\put (29,14) {\small$S_K^\ast$}
	\end{overpic}
	\caption{Region $S_K^*$ (in green) obtained as the intersection of the set of deviations $\mathbb{B}_{\infty}(x^*, \rho)$ around the equilibrium $x^*$ (red dot) with the domain $\Pi_K$. \GP{The samples producing the constraints in blue are in the compression set of $x^*$, while those associated with the red constraint are in the compression set of $S_K^*$; discarding these does not change $x^\ast$.}} \label{final_affine_M=2} 
	\label{fig:setM_3D}
\end{figure}

This is depicted in Figure \ref{final_affine_M=2} using the ${\infty}$-norm (any other norm could have been used instead): an algorithm $\algo$ (see Sec.~\ref{sec:background_scenario_approach}) takes as input a multi-sample $\msample$ and returns the region $S_K^*$ around the solution $x^*\in\mathbb{R}^2$ of a game with two players whose decisions are defined as scalar quantities. \FF{For this pictorial example, $\Pi_K$ is shaped exclusively by sampled coupling constraints. Any compression set as per Definition \ref{def_compression} for $\algo$ must be associated with the solid blue constraints (these form a compression for $x^*$), and with the dashed red constraint that intersects $\mathbb{B}(x^*,\rho)$---as its removal would change $S_K^*$.} 

 We can quantify the number of samples that form a compression set for the algorithm that returns $S_K^*$ in an \emph{a posteriori} fashion as established in Theorem~\ref{thm_support}. To this end, for a fixed confidence $\beta \in (0,1)$, let the violation level be defined as a function $\epsilon: \{0, \dots, K\} \rightarrow [0,1]$ satisfying \cite[Eq.~(7)]{CampiGarattiRamponi2018}
 \begin{equation}
 	\epsilon(K)=1 \quad \text{and} \quad \sum_{i=0}^{K-1} {\binom{K}{i}}(1-\epsilon)^{K-i}=\beta. \label{eps_def}
 \end{equation}

\begin{thm} \label{thm_support}
Under Assumptions \ref{assum_cvx_differ}--\ref{assum_monotonicity_continuity}, let algorithm $\algo$ return a pair $(x^*,S^*_K(x^*))$. Fix  a confidence parameter $\beta \in (0,1)$
and a violation level $\epsilon: \{0, \dots, K\} \rightarrow [0,1]$ that satisfies (\ref{eps_def}). We have that
	\begin{equation}
			\mathbb{P}^{K} \left\{\msample \in \Delta^{K}:~\mathbb{V}(S_K^*)> \epsilon(s^*+M) \right\} \leq \beta, \nonumber
	\end{equation}
where $s^*$ is the number of samples that form a compression set for the equilibrium mapping $x^*$ and $M$ the number of uncertain constraints of $\Pi_K$ that intersect $S_K^*$.
\end{thm}
\emph{Proof}: 
Let $(x^*,S^*_K)$ be the solution returned by $\algo$ for some given $\msample$, according to Definition \ref{Algorithm_definition}. We aim at determining a compression set for $\algo(\msample)$, \FF{and use its cardinality to reach} the theorem's conclusion by means of Theorem 2 in \cite{Pantazis2020}. 
This would be the union of:
(i) the samples that form a compression set for $x^*$---i.e., solving the problem using only these would result in the same equilibrium obtained by using all samples---, and (ii) 
any other sample (not in the compression set of $x^*$) whose removal can still lead to a change of the region $S_K^*$.\\ 
Case (i): Determining a (possibly non-minimal) compression set for $x^*$ can be achieved, as suggested in \cite{CampiGarattiRamponi2018}, by progressively removing samples till a subset that leaves the solution unchanged is determined. We denote its cardinality by $s^*$. With reference to Fig.~\ref{final_affine_M=2}, this set would be associated with the blue constraints active at $x^*$.
Case (ii): We need to count the samples whose removal does not change $x^*$ but yields a larger region $S^*_K$ (red constraint in Fig.~\ref{final_affine_M=2}). Their number can be upper bounded by the $M$ facets of $\Pi_K$ that intersect $S_K^*$. \\
Hence, the number of samples that form a compression set for $\algo(\msample)$ is bounded by $s^*+M$. Existence of a compression set $I$ with a bound on its cardinality is  sufficient for the application of Theorem 2 in \cite{Pantazis2020}. The fact that for the minimal compression set $|I^*|\leq |I|\leq s^*+M$ always holds leads then to the statement of this theorem. \hfill $\blacksquare$ 

It is important to stress that the application of Theorem \ref{thm_support} is agnostic on the choice of the equilibrium seeking algorithm.
To use the result of Theorem \ref{thm_support}, one needs to quantify (an upper bound of) the number of samples $s^*$ that form a compression set for $x^*$ and (an upper bound of) the number $M$ of coupling constraints that correspond to facets of $S_K^*$. 
While $s^*\leq nN$ under Assumptions \ref{assum_cvx_differ}--\ref{assum_monotonicity_continuity},\footnote{By arguments similar to those in \cite{Pantazis2020}, \cite{SchildbachEtAl2012_SupportRank} it can be shown that a tighter bound $s^*\leq n$ holds for the game $G_K$ in case coupling constraints only concern the aggregate variable.} an upper bound for $M$ can in general only be achieved \emph{a posteriori}, i.e., once $\msample$ is sampled. In the next section we show how we could obtain \emph{a priori} bounds for the same quantity.

\subsection{A priori probabilistic certificates}\label{sec_apriori}
Consider the scenario game $G_K$ and suppose that bounded deviations from the solution are allowed. We model such deviations as a ball of radius $\rho$ around the equilibrium, as in Section~\ref{sec_aposteriori}. In contrast to the a posteriori nature of the result therein, our goal here is to achieve an a priori bound. Namely, we aim at establishing the main statement of Theorem 1 with a \emph{prespecified} violation level, which does not depend on the given multi-sample $\msample$.  In other words, we seek a statement---holding with known confidence---of the form $\mathbb{V}(\Pi_K \cap \mathbb{B}(x^*, \rho)) < \bar{\epsilon}$, with $\bar{\epsilon} \in (0,1)$ a priori fixed.
	
To achieve this we build upon the previous conclusions, which expose a link between the probability of constraint violation and the number $M$ of facets of $\Pi_K$ (each originated from some uncertainty sample) that $\mathbb{B}(x^*, \rho)$ intersects.
In particular, a monotonic relationship follows from \eqref{eps_def}: the smaller $M$ the better, i.e., less conservative, the theoretical feasibility guarantees on constraint violation for the strategies belonging to the feasible region $S_K^*$ surrounding the equilibrium. Also, a smaller value of $M$ can result in a larger region for which the guarantees of Theorem 1 hold---due to a smaller portion of $\mathbb{B}(x^*, \rho)$ being cut off by intersection with $\Pi_K$. This motivates us to study the role of $M$ as a modulating parameter for the robustness of the feasibility certificates offered for the region $S_K^*$, as well as the extent of deviation from the nominal equilibrium covered by such certificates.

\subsubsection{GNE-seeking algorithm with a priori robustness guarantees}
\label{sec:GNEapriori}
We consider an iterative scheme to determine a solution of VI$_K$ in \eqref{eq:VI_K}. In particular, since
\FF{the problem involves coupling constraints, we build our Algorithm~\ref{Algorithm_apriori} upon a primal-dual scheme, where constraint satisfaction is achieved by the use of Lagrange multipliers; similar developments hold for both GNE and GWE problems}. To this end, we define the augmented \FF{vector $y:=(x, \mu) \in \mathbb{R}^{nN+m}$ by stacking the} global decision vector $x$ and the Lagrange multipliers $\mu=(\mu_{\ell})_{\ell=1}^m \in \mathcal{M} \subseteq \mathbb{R}_+^m$.  The set $\mathcal{M}$ denotes the domain of $\mu$; in the sequel we impose some structure on $\mathcal{M}$ once some necessary theoretical ingredients are introduced.
As deterministic constraints do not play a role in the evaluation of the robustness guarantees, suppose for ease of exposition that $\Pi_K$ only comprises uncertain coupling constraints. Let $A\in \mathbb{R}^{m \times nN}$ and $b\in\mathbb{R}^{m}$ such that
\begin{subequations}\label{eq:H_repr}
	\begin{gather}
	\Pi_K=\{x \in X: A x \leq b\},\label{eq:H_repr_a} \\
	\|a_{\ell}\|_2 = 1,\; \ell=1,\ldots,m, \label{eq:H_repr_b}
	\end{gather}
\end{subequations}
where $a_{\ell}\trasp$ denotes the $\ell$-th row of $A$. Eq.~\eqref{eq:H_repr} is the irredundant $H$-representation of the polytopic feasibility region $\Pi_K$ defined in \eqref{eq:VI_K}, where the rows of matrix $A$ are unit vectors. Property \eqref{eq:H_repr_b} is key to the second statement in Lemma 2. It entails no loss of generality, since for any $A,b$ forming an equivalent $H$-representation of $\Pi_K$, \eqref{eq:H_repr} can be obtained by normalising each row of $A$ and the corresponding component of $b$ by the row-vector norm.
Thus, the pair $(A,b)$ encodes the set of randomised coupling constraints that constitute facets of $\Pi_K$\footnote{Formally, $A: \Delta^K \rightarrow \mathbb{R}^{m \times nN}$ and $b: \Delta^K \rightarrow \mathbb{R}^{m}$ are mappings from the $K$-multisample to the space of real $m \times nN$ matrices and $m$-dimensional vectors, respectively.}.
\begin{algorithm}[t] 
	\caption{A priori robust GNE seeking algorithm}
	\label{Algorithm_apriori}
	\begin{algorithmic}[1]
		\Require $ y^{(0)},\rho\in\mathbb{R}_{+}, \msample\in\Delta^K, M\leq m, \xi\geq 0$
		\State $\kappa \leftarrow 0$
		\Repeat 
		\State $y^{(\kappa+1)} \leftarrow \text{proj}_{X \times \mathcal{M}, D} \left [ y^{(\kappa)}-D^{-1}T(y^{(\kappa)}, \rho, M) \right ]$\label{alg:update_dyn}
		\State $\kappa \leftarrow \kappa+1$
		\Until $\|y^{(\kappa+1)}-y^{(\kappa)}\|\leq \xi$
		\State $y^* \leftarrow y^{(\kappa+1)}$\\
		\Return $y^*$ and $\Pi_K \cap \mathbb{B}(x^*, \rho)$
	\end{algorithmic}
\end{algorithm} 

\FF{The main step of Algorithm~\ref{Algorithm_apriori} (step 3) is a projected gradient descent (ascent) update for $x$ ($\mu$) through the mapping $T: \mathbb{R}^{nN+m+1} \times \mathbb{N} \rightarrow \mathbb{R}^{nN+m}$ given by}
 \begin{equation}\label{eq:Tmap}
	T(y,\rho,M)\coloneqq\begin{bmatrix}
		F(x)+A\trasp\mu \\
		-\left(Ax-b+Q(\mu, M)\vtight\right)
	\end{bmatrix}.
\end{equation}
\FF{$T$ follows from the primal-dual conditions of the game solution; see \cite[Sec.~4.2]{FacchineiKanzow2009}, \cite[Sec.~1.4.1]{Pang1}. $F$ is the pseudo-gradient mapping defined as in Section~\ref{sec_gamesolVI}, $A,b$ are as in \eqref{eq:H_repr}, and $\vtight\coloneqq\tight\mathbf{1}_m$, where $c$ is a constant scaling factor (see Sec.~\ref{sec_mapQ}) and $M$ a nonnegative integer.
In the second block-row of \eqref{eq:Tmap}, the $m- M$ least relevant (based on the multipliers value) coupling constraints are tightened by an amount $c\rho$ through the mapping $Q: \mathbb{R}^m_+ \times \mathbb{N} \rightarrow \{0,1\}^{m\times m}$. 
Finally, the \emph{asymmetric projection matrix} $D \succ 0$ includes the step-size parameter $\tau>0$ and is defined as
	\begin{equation}\label{eq:D}
	D\coloneqq\begin{bmatrix}
	\frac{1}{\tau}I_{nN}& 0\\
	-2A &  \frac{1}{\tau}I_{m}
	\end{bmatrix}.
	\end{equation}
%
Note that the constraint tightening performed in the second block-row of $T$ is equivalent to preventing $\mathbb{B}(x^*,\rho)$ from intersecting these constraints.  
In other words, $Q$ ensures that the number of facets of $\Pi_K$ intersecting $\mathbb{B}(x^*,\rho)$ is at most $M$, which in turn enables to obtain an \emph{a priori} estimate of the number of samples that form a compression for $S_K^*$ and hence on $\mathbb{V}(\Pi_K \cap \mathbb{B}(x^*, \rho))$; this is formalised by Theorems \ref{thm_convergence} and \ref{thm:apriori_guarantees}.  
Since $m-M$ coupling constraints are tightened, smaller values for $M$ can result in a more robust and possibly larger region $S_K^*$; however, they can also move the location of the nominal equilibrium $x^*$ to a somewhat less efficient point towards the interior of $\Pi_K$. As we will demonstrate numerically in the sequel, this is the case with \emph{potential} games  \cite{FacchineiEtAl2011_Potential}.}

\subsubsection{Constraint tightening via mapping $Q$}
\label{sec_mapQ}
We define the mapping $Q$ as 
\begin{equation}\label{eq:Q}
	Q(\mu,M)\coloneqq P\trasp(\mu)R(M), 
\end{equation} 
where
\begin{itemize}
	\item $P: \mathbb{R}^m \rightarrow \{0,1\}^{m \times m}$ returns a permutation matrix such that $P(\mu)\mu$ is the vector composed by the elements of $\mu$ arranged in decreasing order.
	\item $R: \mathbb{N} \rightarrow \{0,1\}^{m \times m}$ takes as input the number of coupling constraints $M\leq m$ we allow $\mathbb{B}(x^*, \rho)$ to intersect with and returns as output the matrix
	\begin{equation}\label{eq:Rmap}
	R(M)=\begin{bmatrix}
	\begin{matrix}
	\bigzero_{m \times M}
	\end{matrix}
	& \rvline & \begin{matrix}
	0_{M \times m-M}\\
	I_{m-M}
	\end{matrix}
	\end{bmatrix}. 
	\end{equation}
	Compatibly with the definition of $P(\cdot)$, $R(M)P(\mu)\vtight = (\mathbf{0}_M\trasp \; \tight\mathbf{1}_{m-M}\trasp)\trasp = R(M)\vtight$, where the last equality holds since all components of $\vtight$ are equal.
\end{itemize} 
\FF{As discussed in Section~\ref{sec:GNEapriori}, $Q(\cdot,M)$ allows to tighten the constraints corresponding to the smallest $m-M$ multipliers. For this, we use the radius of the sphere that circumscribes $\ball$. This is $\vtight_{\ell} = c\rho\|a_{\ell}\|_2 = c\rho$, where the last equality is due to \eqref{eq:H_repr_b}; depending on the choice of norm,} $c=1$ if $\ball[\cdot]$ is expressed by a $p$-norm with $p\leq 2$, $c=\sqrt{n}$ otherwise. Conversely, at most $M$ constraints can intersect $\mathbb{B}(x^*, \rho)$ upon convergence of the algorithm. Let $\mathcal{L}(M)\subseteq\{1,\ldots,m\}$ contain the indices of the $M$ largest multipliers. Then, $\ell \in \mathcal{L}(M) \Leftrightarrow (Q(\mu,M)\vtight)_\ell = 0$, and the second block row of $T$ in \eqref{eq:Tmap} expresses
\begin{equation}\label{eq:constr_tight}
	\begin{cases}
	a_{\ell}\trasp x \leq b_{\ell} & \text{if $(Q(\mu,M)\vtight)_{\ell}=0$},\\
	a_{\ell}\trasp x \leq b_{\ell} - \tight & \text{if $(Q(\mu,M)\vtight)_{\ell}=c\rho$}.
	\end{cases}
\end{equation}
%
\paragraph*{Illustrative example:}
\FF{Let $\Pi_K$ result from the intersection of 3 hyperplanes and allow $\ball[\cdot]$ to intersect at most $M = 1$ of them. From \eqref{eq:Rmap},  $R(M)=\left[\begin{smallmatrix}
0 & 0 & 0 \\
0 & 1 & 0 \\
0 & 0 & 1	
\end{smallmatrix}\right]$. 
At iteration $\kappa$ of Algorithm \ref{Algorithm_apriori}, let the multiplier vector $\mu^{(\kappa)}=(\mu^{(\kappa)}_{\ell})_{\ell=1}^3$ be such that $\mu^{(\kappa)}_{2}>\mu^{(\kappa)}_{1}>\mu^{(\kappa)}_{3}$.}
Then, $P(\mu^{(\kappa)})= \left[\begin{smallmatrix}
		0 & 1 & 0\\
		1 & 0 & 0 \\
		0 & 0 & 1	
	\end{smallmatrix}\right]$
\FF{is the permutation matrix such that $P(\mu^{(\kappa)}) \mu^{(\kappa)} = (\mu^{(\kappa)}_{2}\; \mu^{(\kappa)}_{1}\; \mu^{(\kappa)}_{3})\trasp$.
So $Q(\mukappa,M)\vtight=P\trasp(\mukappa)R(M)\vtight=(\tight \; 0 \; \tight)\trasp$, where $P\trasp(\cdot)$ applies the correct ordering to the vector $R(M)\vtight$. 
Suppose $\mu^{(j)}_{2}>\mu^{(j)}_{1}>\mu^{(j)}_{3}$ holds for all $j\geq\kappa$. Then, at convergence, it follows from \eqref{eq:constr_tight} that $\ball[x^*]$ will \emph{not} be intersecting the constraints associated to $\mu_1$ and $\mu_3$, whereas it could be intersecting the hyperplane associated to $\mu_2$.}
 
\subsection{Convergence analysis and main result}
\label{sec:mainresults}
\KM{Due to $Q$, mapping $T$ is discontinuous on $X \times \mathbb{R}^m$. To circumvent this, we restrict the multipliers to the set $\mathcal{M}$ on which we impose some structure granting continuity of $T$ on $X \times \mathcal{M}$.
To this end, let $\mathcal{Z}:= [\zeta,+\infty) \cup \{0\}$, for some small $\zeta>0$, i.e., $\mathcal{Z}\subset\mathbb{R}$ contains all nonnegative scalars which take value greater than $\zeta$ when nonzero. 
\begin{assumption} \label{ass:ass_setM}
Let $\Lambda$ be an arbitrarily large compact set. $\mathcal{M}$ admits the form
\begin{multline}\label{eq:M}
	\mathcal{M}:=\left\{ \mu \in \Lambda: (P(\mu)\mu)_{\ell+1} \leq (P(\mu)\mu)_{\ell}-\zeta,\right.\\
				 \left.\forall \ell=1,\dots, m-1\right\} \cap \mathcal{Z}^m.
\end{multline}
\end{assumption}}
Recalling that $P(\mu)\mu$ rearranges the multipliers in descending order, the set $\mathcal{M}$ contains all vectors where the difference between every pair of strictly positive components---and the distance of the smallest of these from zero---is no less than $\zeta$. We note that \eqref{eq:M} \FF{is the union of $q=m!+m+1$ disjoint convex subsets of $\mathbb{R}^m_+$, each of which we denote as $\setM_j$, i.e., $\setM = \bigcup_{j=1}^{q} \setM_j$; figure~\ref{fig:setM_3D} illustrates this set for $m=3$. It is therefore possible to compute the projection in line 3 of Algorithm~\ref{Algorithm_apriori} by, e.g., projecting on $\setM_j$, for $j=1,\ldots,q$, and then setting $y^{(\kappa+1)}$ to be the solution among these that results in the minimum distance from $y^{(\kappa)}-D^{-1}T(y^{(\kappa)}, \rho, M)$. Still, the projection on $\mathcal{M}$ can be computationally intensive if $q$ is large.} 

\begin{figure}
	\centering
	\includegraphics[width=0.9\columnwidth,trim=0 0.6cm 0 0.6cm,clip]{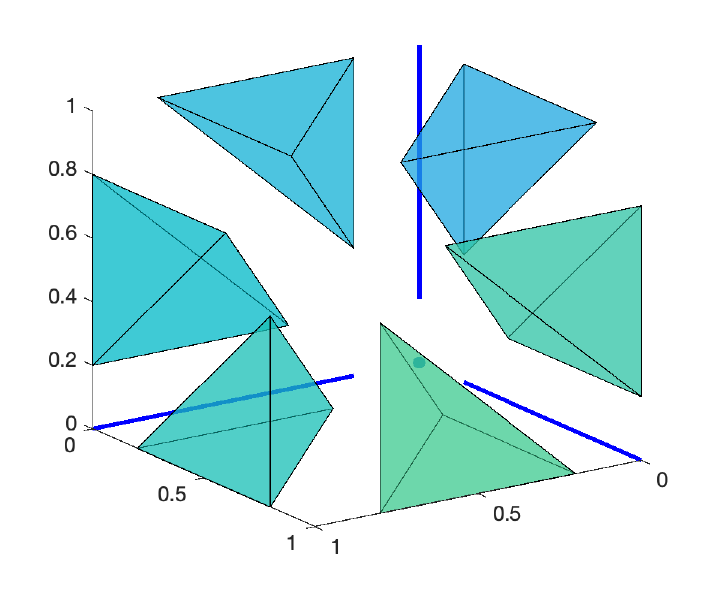}
	\caption{Domain $\setM$ of the Lagrange multipliers associated to the coupling constraints, for  $\zeta = 0.2$ and $m=3$. This results in $q=10$ convex subsets, including the origin and a portion of the axes.} 
	\label{fig:setM_3D}
\end{figure}

Imposing on $\mathcal{M}$ the structure of \eqref{eq:M} endows $T$ with the desired nonexpansiveness properties that are exploited in the proof of Lemma~\ref{lem:Tproperties}. In the numerical implementation of the algorithm, ensuring $\mu \in \mathcal{M}$ can possibly introduce small perturbations in the multipliers---compared to standard formulations where $\mu\in\mathbb{R}^m_+$---which in turn could produce a slight violation of the constraints (this can be controlled through the magnitude of $\zeta$). \FF{We note that $\setM$ is compact by construction due to the intersection with the compact set $\Lambda$ in \eqref{eq:M} which can, however, be arbitrarily large thus not impacting the result numerically. Compactness is used in the proof of Theorem~\ref{thm_convergence}; Remark \ref{rem:compact_relax} discusses cases where this requirement can be lifted.}
 
\begin{lemma}\label{lem:Qbound}
Define $T$ as in \eqref{eq:Tmap}--\eqref{eq:Rmap}, where $A,b$ satisfy \eqref{eq:H_repr}. Then, for any $\mu, \mu' \in \mathcal{M}$, $\mu\neq \mu'$, there exists an integer $0\leq h \leq M$ such that 
	\begin{equation}
	(\mu-\mu')\trasp(Q(\mu,M)-Q(\mu',M))\vtight \leq - h\zeta c\rho.
	\end{equation}
\end{lemma}
\begin{lemma}\label{lem:Tproperties}
	Consider $T$ as in \eqref{eq:Tmap}--\eqref{eq:Rmap}, where $A,b$ satisfy \eqref{eq:H_repr} and $\mathcal{M} = \bigcup_{j=1}^{q} \mathcal{M}_j$ as in \eqref{eq:M}. For each $j=1,\ldots,q$, let VI$(\setX\times\setM_j,T)$ denote the VI problem defined by the map $T$ restricted to the subdomain $\setX\times\setM_j$. 
	Under Assumptions~\ref{assum_cvx_differ}--\ref{assum_monotonicity_continuity} the following holds:
	\begin{enumerate}
	\item $T$ is continuous on $\setX \times \mathcal{M}$.  
	\item Let $D$ as in \eqref{eq:D}
	and set $\tau>0$ such that
	\begin{multline}\label{eq:tau_cond}
	\tau < \min \bigg \{ \frac{-L_F^2 \! +  \! \sqrt{L_F^4+4\alpha^2\|A\|^2}}{2\alpha\|A\|^2}, \\
	\frac{-\rho(1\!+\! \|A\|^2) \! + \! \sqrt{\rho^2(1\!+\! \|A\|^2)^2+16\zeta^2 \|A\|^2}}{4\zeta\|A\|^2 } \bigg \}.
	\end{multline}
	Then, for any $j=1,\ldots,q$, Algorithm~\ref{Algorithm_apriori} converges to a solution of VI$(\setX\times\setM_j,T)$, when the gradient step in line \ref{alg:update_dyn} is projected on the corresponding subdomain, for any $y^{(0)}\in\setX\times\setM_j$.
	\end{enumerate}
\end{lemma}
Continuity of the mapping is essential for the theoretical convergence of Algorithm 1. The second part of Lemma~\ref{lem:Tproperties} provides an admissible range of values for $\tau$ such that Algorithm~\ref{Algorithm_apriori}  converges to a solution of VI$(\setX\times\setM_j,T)$ if at each iteration the projection in line \ref{alg:update_dyn}  is performed on the (convex) subdomain $\setM_j \subset \setM$, $j\in\{1,\ldots,q\}$. \FF{The stepsize $\tau$ is chosen such that conditions standard in NE seeking are satisfied and oscillations among multiple equilibria are avoided.}
Still, we are interested in establishing convergence on the entire domain  $\setM$, so at each iteration the projected solution might belong to a different subdomain. This does not trivially follow from the second part of Lemma~\ref{lem:Tproperties}; therefore, by Lemmas~\ref{lem:Qbound} and~\ref{lem:Tproperties} we establish an additional condition on $\tau$ such that Algorithm~\ref{Algorithm_apriori} retrieves a solution of VI$(\setX\times\setM,T)$. 

\begin{thm}\label{thm_convergence}
Consider Assumptions \ref{assum_cvx_differ}, \ref{assum_affine}, \ref{assum_monotonicity_continuity} and \ref{ass:ass_setM}. Fixed $0\leq M\leq m$, assume the domain $\Pi_K$ is nonempty for any of the ${m \choose M}$ combinations of constraints tightened as in \eqref{eq:constr_tight}. Let $D\succ 0$ be defined as in \eqref{eq:D}, where $\tau$ satisfies \eqref{eq:tau_cond} and
	\begin{equation}\label{eq:thm_bound_tau}
	\tau <  \frac{-(\bar{C}+\bar{R})+ \sqrt{\vphantom{\big(}(\bar{C}+\bar{R})^2 +2\zeta\bar{R}}}{2\bar{R}},
	\end{equation}
where $\bar{R}=\max \big \{\sup_{x\in\setX} \sup_{\mu \in \setM}\|2A(F(x) +A\trasp \mu)\|, \\ \sup_{x\in\setX}\|Ax-b\|\big \}$ and $\bar{C} = \tight \sqrt{m-M}$. 

Then Algorithm \ref{Algorithm_apriori} converges to a solution of VI$(\setX\times\setM)$ for any initial condition $y^{(0)}\in\setX\times\setM$.
\end{thm}
Note that as $\mukappa\to\mu^*$, we have $Q(\mukappa)\to Q(\mu^*) =: Q^*$. Then, the solution returned by Algorithm 1 is the equilibrium of a variant of $G_K$ with $m-M$ tightened constraints (follows from \eqref{eq:constr_tight} with $Q(\mu)$ replaced by $Q^*$). 
\begin{remark}[Relaxing compactness]\label{rem:compact_relax}
\KM{Theorem~\ref{thm_convergence} still holds when $\Lambda=\mathbb{R}^m$ in the definition of  $\setM$ in \eqref{eq:M} if for all multi-samples, $(i)$ $A$ is full row-rank, or $(ii)$ all elements of $A$ are positive. \\
$(i)$
To show this, consider mapping $T$ and matrix $D$ in \eqref{eq:D}. The multipliers' update involves projecting (weighted according to $D$) on $\mathcal{M}$, the term
\begin{multline}\label{eq:integr}
 \mu^{(k)} - 2\tau^2 A(F(x^{(k)}) +A\trasp \mu^{(k)})  \\
 +\tau (Ax^{(k)} - b + Q(\mu^{(k)},M) \vtight). 
\end{multline}
Since $X$ is compact, there exists a subsequence $\{\kappa_i\}_{i \in \mathbb{N}}$ such that $\lim_{i\to \infty} \kappa_i = \infty$, $\lim_{i\to \infty} x^{(\kappa_i)} = \overline{x}$, for some $\overline{x} \in X$.
It suffices to show that the sequence of multipliers 
$\{\mu^{(\kappa_i)}\}_{i=1}^{\infty}$ remains bounded (all arguments in the proof of Theorem \ref{thm_convergence} from \eqref{eq:wlimitset} onwards remain unaltered). \\
For the sake of contradiction, assume that there exists at least one element of $\mu^{(\kappa_i)}$ that tends to infinity across the considered subsequence. 
Let then $\mu^{(\kappa_i)} = (\mu_{\infty}^{(\kappa_i)}~ \mu_{F}^{(\kappa_i)})$, where based on our contradiction hypothesis $\lim_{i\to \infty} \| \mu_{\infty}^{(\kappa_i)} \| = \infty$ while $\lim_{i\to \infty} \| \mu_{F}^{(\kappa_i)} \| < \infty$. (Taking the first elements of $\mu^{(\kappa_i)}$ to be the ones that tend to infinity is only to simplify notation and is without loss of generality.) Let then $A = [A_{\infty} \trasp ~ A_{F} \trasp] \trasp$, $b = (b_{\infty}~ b_{F})$ be the corresponding partition of $A$ and $b$, respectively, where $A_{\infty}, b_{\infty}$ are non-empty by hypothesis.
To have $\|\mu_{\infty}^{(\kappa_i)}\|\to\infty$, we need the terms that are integrated in the multipliers' update, i.e., last two terms in \eqref{eq:integr}, to be positive for all $i$ (in fact across a subsequence), which since $\tau>0$ is equivalent to
\begin{multline}\label{eq:integr1}
A_{\infty}  A\trasp \mu^{(\kappa_i)} < - A_{\infty}  F(x^{(\kappa_i)})  \\
+ \frac{1}{2\tau} \Big (A_{\infty}  x^{(\kappa_i)} - b_{\infty} + (Q(\mu^{(\kappa_i)},M) \vtight)_{\infty} \Big ),
\end{multline}
where $( \cdot )_{\infty}$ denotes the elements of its argument corresponding to $\mu_{\infty}^{(\kappa_i)}$. Notice that $A_{\infty}  A\trasp \mu^{(\kappa_i)}  = A_{\infty}  A_{\infty}\trasp \mu_{\infty}^{(\kappa_i)} + A_{\infty} A_{F}\trasp \mu_{F}^{(\kappa_i)}$. As such, we have
\begin{multline}\label{eq:integr2}
A_{\infty} A_{\infty}\trasp \mu_{\infty}^{(\kappa_i)} < - A_{\infty}  A_{F}\trasp \mu_{F}^{(\kappa_i)} - A_{\infty}  F(x^{(\kappa_i)})  \\
+ \frac{1}{2\tau} \Big (A_{\infty} x^{(\kappa_i)} - b_{\infty} + (Q(\mu^{(\kappa_i)},M) \vtight)_{\infty} \Big ).
\end{multline}
However, $\lim_{i\to \infty} x^{(\kappa_i)} = \overline{x} \in X$, and $(Q(\mu^{(\kappa_i)},M) \vtight )_{\infty} \leq c \rho$ for all $i$, while by Lemma \ref{lem:Tproperties}, $F$ is continuous over the domain of multipliers satisfying \eqref{eq:M}. Moreover, $\mu_{F}^{(\kappa_i)}$ contains the components of $\mu^{(\kappa_i)}$ that remain finite. Therefore, the limit as $i\to \infty$ of the right-hand side of \eqref{eq:integr2} is finite. Due to the assumed full row-rank structure of $A$, matrix $A_{\infty} A_{\infty}\trasp$ is invertible, hence \eqref{eq:integr2} implies 
$\lim_{i\to \infty} \| \mu_{\infty}^{(\kappa_i)} \| < \infty$, establishing a contradiction showing that the subsequence $\{\mu^{(\kappa_i)}\}$ remains bounded. \\$(ii)$ If all elements of $A$ are positive, and since $a_{\ell} \trasp a_\ell = 1$, for all $\ell =1,\ldots,m$, all arguments of case (i) remain the same with the only difference that we directly have that $\|A_{\infty} A_{\infty}\trasp \mu_{\infty}^{(\kappa_i)}\| \geq \|\mu_{\infty}^{(\kappa_i)}\|$.}
\end{remark}

The next result accompanies the region $S_K^*= \Pi_K \cap \ball$ of strategic deviations from the equilibrium $x^*$ with \emph{a priori} probabilistic feasibility guarantees that can be tuned by means of $M$. 
It should be noted that Theorem \ref{thm_convergence} establishes that there exists a choice of $\tau$ to guarantee convergence of Algorithm \ref{Algorithm_apriori}. The admissible range of values for $\tau$ is explicit via  \eqref{eq:tau_cond}, \eqref{eq:thm_bound_tau}, but difficult to quantify due to $\bar{R}$. Numerical evidence suggests that selecting a small enough value is sufficient for convergence.

\begin{thm} \label{thm:apriori_guarantees}
Consider Assumptions \ref{assum_cvx_differ}, \ref{assum_affine}, \ref{assum_monotonicity_continuity} and \ref{ass:ass_setM}. 
Let $x^*$ and $S_K^* = \Pi_K \cap \ball$ be returned by Algorithm~\ref{Algorithm_apriori}; fix $\overline{\epsilon} \in (0,1)$ and $M$. We then have that
	\begin{multline}\label{eq:Theorem3_guarantees}
		\mathbb{P}^{K} \Big \{\msample \in \Delta^{K}:\; \mathbb{V}(S_K^*)\leq \overline{\epsilon}\Big \}\\
		 \geq 1-\sum_{i=0}^{nN+M-1} \binom{K}{i}\overline{\epsilon}^i(1-\overline{\epsilon})^{K-i}. 
	\end{multline}
\end{thm}

By Definition 2, Theorem~\ref{thm:apriori_guarantees} guarantees that for any point in $S_K^*$, the probability of constraint violation is bounded by $\bar{\epsilon}$, with confidence at least $1-\sum_{i=0}^{nN+M-1} \binom{K}{i}\overline{\epsilon}^i(1-\overline{\epsilon})^{K-i}$. The dependence of this term on $M$ gives us an additional degree of freedom in trading the robustness of the solution for its associated probabilistic confidence. The choice of $M$ can also have an effect on the size of $S_K^*$, as well as on the location of $x^*$, thus resulting in a trade-off between performance and robustness.

For the case in which the coupling constraints concern exclusively the aggregate variable, it can be shown that the upper limit of the summation in the right-hand side of \eqref{eq:Theorem3_guarantees} can be replaced by $n+M-1$, as $n$ is the dimension of the aggregate vector. This allows to state \eqref{eq:Theorem3_guarantees} with a much higher confidence of $1-\sum_{i=0}^{n+M-1} \binom{K}{i}\overline{\epsilon}^i(1-\overline{\epsilon})^{K-i}$; for details, we refer the reader to \cite{Pantazis2020}, where the notion of support rank is exploited \cite{SchildbachEtAl2012_SupportRank}.

\section{Numerical example}\label{sec:num_example}
Consider a game with $N$ agents whose decisions are subject to deterministic local constraints and uncertain coupling  constraints on the aggregate decision: 
\begin{equation}\label{eq:example_scenario}
	\left.
	\begin{aligned}
		& \min_{x_i \in X_i} \ x_i\trasp (C\sigma(x)+d) \\  
		& \text{ subject to } \  \lowb \leq \sigma(x) \leq \uppb, \\ 
		& \qquad \qquad \qquad \qquad \qquad k=1, \dots, K 
	\end{aligned}
	\right\}\; \forall i\in\mathcal{N},
\end{equation}
where $C \succ \alpha I_n$, for some $\alpha>0$, and $d \in \mathbb{R}^n$. Note that a structure similar to our numerical example has been considered  in applications of aggregative games such as  electric vehicle charging  and traffic management under uncertainty \cite{Dario2019}, \cite{Feleconf2019}, \cite{Fele2021}. We impose no knowledge of $\Delta$ and $\mathbb{P}$; we rely instead on a scenario-based approximation of the game, whereby each sample $\delta_k\in\msample$ gives rise to $\lowb,\uppb$. Eq.~\eqref{eq:example_scenario} is an \emph{aggregative} game in the form of \eqref{eq:scenario_game}. 
In this instance, we assume each agent's action has negligible effect on the aggregate, and accordingly consider a GWE-seeking problem. Following the definition of $F_{\mathrm{WE}}$ (Sec.~\ref{sec_gamesolVI}), we get $F(x)=F_{\mathrm{WE}}=[C\sigma(x)+c]_{i \in \mathcal{N}}$. It can be verified that $F$ is Lipschitz continuous and strongly monotone with respect to  $\sigma$: by \cite[Thm.~2.3.3]{Pang1}, \eqref{eq:example_scenario} admits a unique aggregate equilibrium $\sigma^*=\sigma(x^*)$.\footnote{We note that this case slightly transcends the conditions in Theorem~\ref{thm_convergence}, as $F$ does not comply with Assumption~\ref{assum_monotonicity_continuity}-\eqref{strong_monotonicity}. Convergence of Algorithm~\ref{Algorithm_apriori} (following from the nonexpansiveness of $T$ on each subdomain $\setM_j$) can still be ensured here due to the affine structure of $F$; cf.~\cite[Sec.~12.5.1]{Pang1}.
}

\begin{figure*} 
	\centering
	\begin{subfigure}[b]{0.33\textwidth} 
		\centering
		\includegraphics[trim=0 0 0 0.5mm, clip,width=\columnwidth]{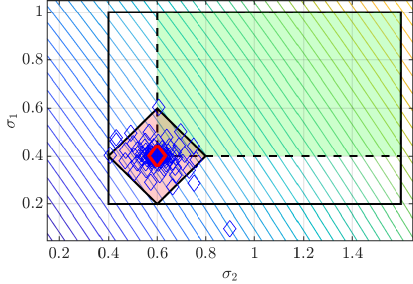}
		\caption{$M=0$} 
		\label{fig:convergence1} 
	\end{subfigure}
	\hfill
	\begin{subfigure}[b]{0.33\textwidth}
		\flushright
		\includegraphics[trim=0 0 0 0.5mm, clip,width=\columnwidth]{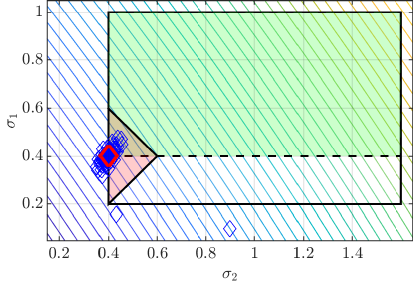} 
		\caption{$M=1$}
		\label{fig:convergence2}
	\end{subfigure}
	\hfill
	\begin{subfigure}[b]{0.33\textwidth}
		\flushleft
		\includegraphics[trim=0 0 0 0.5mm, clip,width=\columnwidth]{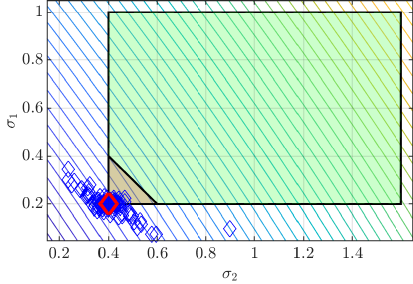}
		\caption{$M\geq 2$} 
		\label{fig:convergence3}
	\end{subfigure}
	\caption{Iterates generated by Algorithm~\ref{Algorithm_apriori}  (blue diamonds) for different choices of $M$. In this numerical instance, $N=50$, $\rho = 10$, and $X_i:=\{x_i \in \mathbb{R}^n: x_i \in [\underline{x}_i, \overline{x}_i] \}$, with $\underline{x}_i=(0 \  0)$, $\overline{x}_i=(3.5 \ 3.5)$. The randomly generated coupling constraints form the rectangular feasibility region $\Pi^\sigma_K$ (delineated by the solid black line). The red-shaded region represents the intersection between the latter and the ball $\mathbb{B}_1(\sigma^*, \rho/N)$ around the aggregate equilibrium $\sigma^*$ (red diamond marker). In this instance, its volume increases as larger values for $M$ are chosen. The value associated to the contour lines of the potential function $E$ decreases from top-right to bottom-left, showing that $\sigma^*$ is the unique minimiser in the admissible region (shaded in green) after constraint tightening is performed by the algorithm (see Sec.~\ref{sec_mapQ}).}
	\label{fig:convergence}
\end{figure*}

We employ Algorithm \ref{Algorithm_apriori} to seek a WE $x^*$ such that, by fixing $M$, a prespecified theoretical violation level is guaranteed for the set $\Pi_K \cap \mathbb{B}(x^*,\rho)$. Due to uniqueness of $\sigma^*$, all sets $\mathbb{B}(\cdot, \rho)$---parametrised by any $x^*$ solving \eqref{eq:example_scenario}---are projected on the unique ball  $\mathbb{B}(\sigma^*, \rho/N)$ in the aggregate space. Also note that by definition of $\sigma$, at most $n$ non-redundant samples will contribute to define the domain $\Pi^{\sigma}_K:=\tfrac{1}{N}(\mathbf{1}_N\trasp\otimes I_n)X \cap \left(\bigcap_{k=1}^{K} C_{\delta_k}\right)$ in \eqref{eq:example_scenario}. For the derivation of the robustness guarantees, we can thus restrict our attention to $S_K^*=\Pi^{\sigma}_K \cap \mathbb{B}(\sigma^*, \rho/N)\subseteq\mathbb{R}^n$. As remarked at the end of Section~\ref{sec:mainresults}, we can apply \eqref{eq:Theorem3_guarantees} with the upper limit in the summation involved replaced by $n+M-1$.
For the case $n=2, N=50$, and \FF{different choices of $M$, Figure~\ref{fig:convergence} depicts the projected iterations $\{\sigma(x^{(\kappa)})\}$, $\kappa=1,2,\ldots$ generated by Algorithm~\ref{Algorithm_apriori} on the space $\Pi^{\sigma}_K$.} It can be observed how the region $S_K^*$ changes as the value of $M$ is modified.

\begin{figure}
	\centering
	\includegraphics[trim=0 0.3cm 0 0.3cm,clip, width=\columnwidth]{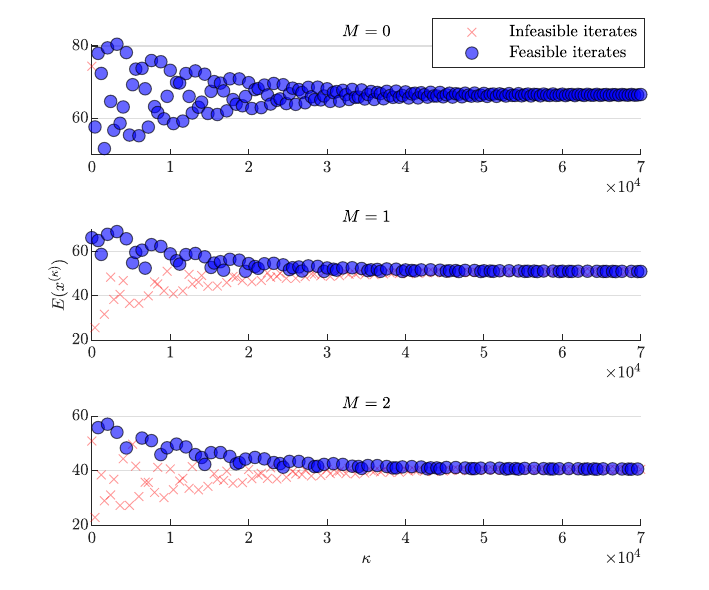}
	\caption{Potential function $E(x^{(\kappa)})$ evaluated along the iterations of Algorithm~\ref{Algorithm_apriori}. Lower values of $M$ yield better confidence on the theoretical robustness certificates for the considered region (see Thm.~\ref{thm:apriori_guarantees}), which results in a lower empirical probability of constraint violation. On the other hand, the system-level efficiency of the equilibrium increases for higher values of $M$. 
	\label{fig:potential_plot}}
\end{figure}

It is worth noting that in this case $F(x)$ is \emph{integrable}---this can be inferred by \cite[Thm.~1.3.1]{Pang1} since the Jacobian of the game is symmetric, i.e., $\nabla_x F(x)=\nabla_x F(x)\trasp$. Therefore, a GWE $x^*$ can also be obtained by solving
\begin{equation}\label{eq:pot_fun}
\begin{aligned}
& \min_{x \in X} \ \sigma(x)\trasp C \sigma(x) + d\trasp \sigma(x) \\
& \text{subject to} \  \lowb \leq \sigma(x) \leq \uppb, \; \  k=1, \dots, K.
\end{aligned}
\end{equation}
In other words, this game admits a \emph{potential function} $E(x):=\sigma(x)\trasp C \sigma(x) + d\trasp \sigma(x)$, whose minimizers correspond to GWEs. $E$ can be interpreted as the total cost incurred by the population of agents, and its minimization leads to the optimum social welfare. The contour lines of $E$ are depicted in Figure~\ref{fig:convergence}: since $x^*$ minimises $E(\cdot)$, $\sigma^*$ lies on the contour associated to the minimum value of $E$ within the feasible domain. Lower values of $M$ result in larger regions for which guarantees are provided. Figure~\ref{fig:potential_plot} shows how the sequence $\{E(x^{(\kappa)})\}_{\kappa=1,2,\ldots}$ converges to the minimum \emph{potential} within the possibly tightened feasibility region. It can be observed how in this case the efficiency of the equilibrium decreases as smaller values of $M$ are chosen. The three panels in Figure~\ref{fig:potential_plot} show the trade-off between system level efficiency and the guaranteed robustness levels. The lower the value of $M$, the lower the empirical constraints violation---corresponding to a better confidence bound in the right-hand side of \eqref{eq:Theorem3_guarantees}.

\section{Concluding remarks}\label{sec:conclusion}
This work proposes a data-driven equilibrium-seeking algorithm such that probabilistic feasibility guarantees are provided for a region surrounding a game equilibrium.  
These guarantees are a priori and the region that is accompanied with such a probabilistic certificate is tunable.
For games that admit a potential function, the proposed scheme is shown to achieve a trade-off between cost and the level of probabilistic feasibility guarantees. In fact, our scheme returns the most efficient equilibrium such that the predefined guarantees are achieved.  Proving this conjecture is left for future work. \KM{Moreover, current work investigates a distributed implementation of the proposed equilibrium seeking algorithm.}

\section{Appendix}
\subsection{Proof of Lemma \ref{lem:Qbound}}
Let $\mu, z$ be arbitrary vectors in $\mathcal{M}$ and, as in the proof of Lemma~\ref{lem:Tproperties}, define $\vec{\mu}, \vec{z}$ as the vectors composed by rearranging the elements of $\mu,z$ in decreasing order. According to this arrangement, let $\mathcal{I}_{\mu}=\{i_1,i_2, \dots, i_m\}$ be the ordered set of indices of $\mu$, i.e., $i_k:\, \mu_{i_k} = \vec{\mu}_{k}$, $k=1,\ldots,m$; as a result, $i_1$ and $i_m$ will be the indices of the largest and smallest components of $\mu$, respectively. Applying a similar definition to $z$, we denote the corresponding set $\mathcal{I}_{z}:=\{j_1,j_2, \dots, j_m\}$. Then, the first $M$ indices in $\mathcal{I}_{\mu}$ and $\mathcal{I}_{z}$, denoted as  $\erre$ and $\erre[z]$, respectively, are relative to the constraints not tightened by the application of $Q(\cdot,M)$. In other words, for all $\ell\in\erre$, $(Q(\mu,M)\vtight)_{\ell} = 0$---and similarly for $z$. Vice versa, the complementary sets $\erre^c = \mathcal{I}_{\mu}\setminus\erre$ and $\erre[z]^c = \mathcal{I}_{z}\setminus\erre[z]$ are such that for all $\ell\in\erre^c$, $(Q(\mu,M)\vtight)_{\ell} = \tight$, and for all $\ell\in\erre[z]^c$, $(Q(z,M)\vtight)_{\ell} = \tight$.
 Let $q=[Q(\mu,M)-Q(z,M)]\vtight$. We distinguish between the following cases:
 \begin{enumerate}
 	\item $\ell \in \erre^c \cap \erre[z]$: we have $(Q(\mu,M)\vtight)_{\ell} = \tight$ since $\ell \in \erre^c$, while $(Q(z,M)\vtight)_{\ell} = 0$ as $\ell \in \erre[z]$. Then, $q_{\ell}=\tight$.
 	\item $\ell \in \erre \cap \erre[z]^c$: from $\ell \in \erre[z]^c$ we have $(Q(z,M)\vtight)_{\ell} = \tight$. On the other hand, since $\ell \in \erre$, $(Q(\mu,M)\vtight)_{\ell} = 0$. This results in $q_{\ell}=-\tight$.
 	\item $\ell \in (\erre \cap \erre[z]) \cup (\erre^c \cap \erre[z]^c)$. If $\ell \in  \erre \cap \erre[z]$ then $(Q(\cdot,M)\vtight)_{\ell} = 0$ for both $\mu$ and $z$. Therefore, $q_{\ell} = 0$. Conversely, if $\ell \in \erre^c \cap \erre[z]^c$, then $(Q(\cdot,M)\vtight)_{\ell} = \tight$ for both $\mu$ and $z$, which results again in $q_{\ell} = 0$.
 	\end{enumerate} 
The sets $\erre^c \cap \erre[z]$, $\erre \cap \erre[z]^c$, $(\erre \cap \erre[z]) \cup (\erre^c \cap \erre[z]^c)$ are pairwise disjoint and exhaust the set $\{1, \dots, m\}$. Hence we can write
\begin{equation}\label{eq:Lemma2_U}
	\begin{split}
	&U=(\mu-z)\trasp(Q(\mu,M)-Q(z,M))\vtight= \sum_{\ell=1}^{m}(\mu_{\ell}-z_{\ell})q_{\ell} \\
	&= \sum_{i \in \erre^c \cap \erre[z]} \!\!\! (\mu_{i}-z_{i})\tight + \!\!\! \sum_{j \in \erre \cap \erre[z]^c} \!\!\! (\mu_{j}-z_{j})\cdot(-\tight)  \\
	&= \bigl(\!\! \underbrace{\sum_{i \in \erre^c\cap \erre[z]} \!\!\! \mu_{i} \; 
	- \!\!\! \sum_{j \in \erre \cap \erre[z]^c} \!\!\! \mu_{j} }_{=:U_1} \; 
	+ \!\!\!  \underbrace{\sum_{j \in \erre \cap \erre[z]^c} \!\!\! z_{j} \;
	- \!\!\! \sum_{i \in \erre^c\cap \erre[z]} \!\!\! z_{i} }_{=:U_2} \bigr)\tight, 
	\end{split}
\end{equation}
Now, notice that for any $i\in\erre^c \cap \erre[z]\subseteq\erre^c$ and $j\in\erre \cap \erre[z]^c\subseteq \erre$, we have by definition of $\erre$ and $\erre^c$ that $\mu_{i} \leq \mu_{j}$ (which by \eqref{eq:M} only holds with equality if $\mu_i=\mu_j=0$). 
With analogous reasoning, we have $z_{i} \geq z_{j}$ for any $i \in  \erre^c \cap \erre[z]\subseteq \erre[z]$ and $j \in \erre \cap \erre[z]^c \subseteq \erre[z]^c$. Let $h_1$ be the cardinality of the set $\erre^c \cap \erre[z]$, and $h_2$ that of $\erre \cap \erre[z]^c$. Then, 
\begin{multline*}
	h_1=|\erre^c \cap \erre[z]|\stackrel{(a)}{=}|\erre[z] \setminus \erre|=|\erre[z]|-|\erre \cap \erre[z]| \\
	\stackrel{(b)}{=}|\erre|-|\erre[z] \cap \erre|=|\erre \setminus \erre[z]|=|\erre \cap \erre[z]^c|=h_2,
\end{multline*}
where $(a)$ holds since $\erre,\erre[z]\subseteq\{1,\ldots,M\}$, and $(b)$ follows from $|\erre|=|\erre[z]|=M$. Therefore $h_1=h_2=:h$ and $0\leq h \leq M$, which implies $U_1\leq 0$ and $U_2\leq 0$ in \eqref{eq:Lemma2_U}. 
We can observe that $U_1<0$ and $U_2<0$ if $\erre \cap \erre[z]^c$ and $\erre^c \cap \erre[z]$ are nonempty and the corresponding components of $\mu$ and $z$ are nonzero. In such a case $h\geq 1$ and we can write
\begin{equation}\label{eq:Lemma2_U1}
	U_1= \sum_{i \in \erre^c\cap \erre[z]} \!\!\! \mu_{i} \; - \!\!\! \sum_{j \in \erre \cap \erre[z]^c} \!\!\! \mu_{j} \leq - h \zeta,
\end{equation}
where the inequality follows from \eqref{eq:M} and the above discussion. A similar reasoning holds for $U_2$. Lastly, note that if $\mu\neq z$ and $h\geq 1$, then at least one of $U_1\leq -h\zeta$ and $U_2\leq -h\zeta$ will hold. By \eqref{eq:Lemma2_U}, we can thus conclude $U \leq -h \zeta \tight$ for any $\mu,z\in\mathcal{M}$, $\mu\neq z$. \hfill $\blacksquare$

\subsection{Proof of Lemma \ref{lem:Tproperties}}
Part (1):  To prove that the mapping $T$ is continuous on its domain, we first notice that $T$ is by construction continuous on $X \times \mathcal{M}$ when the operator $Q(\cdot,M)$ is continuous on $\mathcal{M}$ (as the parameter $M$ is fixed). Therefore, it is sufficient to show that for any $\mu, z \in \mathcal{M}$ and any $\eta > 0$, there exists $\delta>0$ such that 
\begin{equation}\label{eq:lem2_1}
\|\mu-z\| < \delta~ \Rightarrow~ \|Q(\mu,M)-Q(z,M)\|\|\vtight\| < \eta, 
\end{equation}
where $\vtight = \tight\mathbf{1}_{m}\neq \bfzero$.
To this end, consider any $\mu,z \in \mathcal{M}$ such that $\|\mu-z\|< \tfrac{\zeta}{2}$, with $\zeta$ as defined in~\eqref{eq:M}.\footnote{The proof of this part also holds for $\mu\in\mathbb{R}^m_+\supset\mathcal{M}$.} Let $\vec{\mu}$ and $\vec{z}$ denote the vectors $\mu$ and $z$ sorted in decreasing order; thus, $\vec{\mu}_{\ell}$ is the $\ell$-th largest element of $\mu$ (and similarly for $z$). For any given $\ell$, let $i:\, \mu_i = \vec{\mu}_{\ell}$, $j:\, z_j = \vec{z}_{\ell}$, and $\bar{\ell} := \min_{\{1,\ldots,m\}} \ell:\, i \neq j$. In words, $\bar{\ell}$  is the smallest index for which the $\ell$-th largest elements of $\mu$ and $z$ do not appear at the same row of their respective vectors. We then let $\mathcal{I}$ be the set of indices for which the ordering of the elements of $\mu$ and $z$ agrees, i.e., for all $k \in \mathcal{I}$, there exists $\ell<\bar{\ell}$ such that $i=j=k$, with $i:\,\mu_i = \vec{\mu}_{\ell}$ and $j:\,z_j=\vec{z}_{\ell}$.

We prove our statement by contradiction. Suppose there exists $i,j \notin \mathcal{I}$ such that $i:\, \mu_i = \vec{\mu}_{\ell}$ and $j:\, z_j = \vec{z}_{\ell}$ for some $\ell>\bar{\ell}$, where $\mu_i<\mu_j$ and $z_i>z_j$. First, we note that such an instance exists by hypothesis, as otherwise the only possible case is where $i=j$, which contradicts $i,j\notin\mathcal{I}$ and implies $Q(\mu,M)=Q(z,M)$. Since $z \in \mathcal{M}$, it further holds $z_j< z_i-\zeta$, which by $\|\mu_i-z_i\|\leq\|\mu-z\|< \tfrac{\zeta}{2}$ implies
\begin{equation}
	z_j< \mu_i+\frac{\zeta}{2}-\zeta. \label{one}
\end{equation}
We bound \eqref{one} from below by noting $z_j> \mu_j-\tfrac{\zeta}{2}$, which holds since $\|\mu_j-z_j\|<\tfrac{\zeta}{2}$, obtaining $\mu_j-\frac{\zeta}{2} < \mu_i+\frac{\zeta}{2}-\zeta$,
or equivalently $\mu_j<\mu_i$, which contradicts our hypothesis. Hence the elements of any pair of vectors $\mu,z \in \mathcal{M}$ such that $\|\mu-z\|<\tfrac{\zeta}{2}$ must follow the same ordering. By definition of $P(\cdot)$, this implies $P(\mu)=P(z)$ and, in turn, $\|Q(\mu,M)-Q(z,M)\|=0$. This validates \eqref{eq:lem2_1} with $\delta=\tfrac{\zeta}{2}$ and any $\eta>0$, establishing the continuity of $Q(\cdot,M)$ on $\mathcal{M}$ and concluding the proof of the first part.

Part (2): We show that the mapping $T$ fulfils certain nonexpansiveness properties required for the convergence of Algorithm~\ref{Algorithm_apriori}, for compatible choices of $\tau$. In particular, we provide a sufficient condition for which the iteration
\begin{equation}\label{eq:iter_subset}
	y^{(\kappa+1)} = \text{proj}_{\setX \times \setM_j, D} \left [ y^{(\kappa)}-D^{-1}T(y^{(\kappa)}, \rho, M) \right ],
\end{equation}
converges to a solution of VI$(\setX\times\setM_j,T)$, where $j\in\{1,\ldots,q\}$ is fixed, for any $y^{(0)}\in\setX\times\setM_j$. Notice that in \eqref{eq:iter_subset} the skew projection is performed on the convex subdomain $\setX\times\setM_j$.  
\eqref{eq:iter_subset} is the solution of the VI$(\setX\times\setM_j,T_D^{(\kappa)})$ (see \cite[Sec.~12.5.1]{Pang1}), where $T_D^{(\kappa)}(y):= T(y^{(\kappa)},\rho,M) + D(y-y^{(\kappa)})$ is strongly monotone due to $D\succ 0$ and $(T(y,\rho,M)-T(y',\rho,M))\trasp (y-y')\geq 0$, for all $y,y'\in\setX\times\setM$, which in turn follows from Assumption~\ref{assum_monotonicity_continuity} and Lemma~\ref{lem:Qbound}.
The fixed-point iteration \eqref{eq:iter_subset} is an instance of the forward-backward splitting method: we thus resort to standard results in the literature to prove its convergence. Following the notation in \cite[Sec.~12.5.1]{Pang1}, we let $\tilde{D}:= \Dsqr(D-D_s)\Dsqr$, where $D_s:=\frac{D+D\trasp}{2}$. Also, $\setW_j:=\{\Dsqr[] y:\, y\in\setX\times\setM_j\}$, $\setW = \bigcup_{j=1}^q \setW_j$, and $\tilde{T}(w):= \Dsqr T(\Dsqr w, \rho, M)$, for all $w\in\setW$. To ease notation, we drop the dependence of $\tilde{T}$ and $\TD$ on $\rho,M$, as they remain fixed throughout the proof. According to \cite[Thm.~12.5.2]{Pang1} (see also \cite[Sec.~4.3]{ZhuMarcotte1996SIAM}), to ensure convergence of \eqref{eq:iter_subset} to a solution of the VI$(\setX\times\setM_j,T)$ it is sufficient to show that $\TD = \tilde{T}-\tilde{D}$ is $\beta$-cocoercive on $\setW_j$, i.e.,
\begin{equation}\label{eq:cocoerc}
	(\TD(v) - \TD(w))\trasp(v-w)
		\geq \beta\|\TD(v) - \TD(w)\|^2,
\end{equation}
for some $\beta>\tfrac{1}{2}$ and all $v,w\in\setW_j$, $j\in\{1,\ldots,q\}$. In fact, we will go a step forward and demonstrate here that $\TD$ is co-coercive on $\setW$ with $\beta>\tfrac{1}{2}$. Due to the saddle problem structure of the mapping in \eqref{eq:Tmap}, we adopt the procedure in \cite[Prop.~12.5.4]{Pang1} and define $D$ as in \eqref{eq:D} (see also \cite{Dario2019}). It then follows from the above definitions that $\TD(w)$, for any $w\in\setW$, reduces to
\begin{equation}\label{eq:TD}
	\TD(\Dsqr[] y) = \Dsqr\begin{bmatrix}
	F(x) \\
	b-Q(\mu,M)\vtight
	\end{bmatrix}, \; \forall y\in\setX\times\setM,
\end{equation}
which can be easily seen by rewriting \eqref{eq:Tmap} as
\begin{equation*}\label{eq:Tmap_bis}
T(y,\rho,M)=\begin{bmatrix}
F(x) \\
0
\end{bmatrix}+\underbrace{
	\begin{bmatrix}
	0 & A\trasp \\
	-A & 0
	\end{bmatrix}}_{D-D_s}y+\begin{bmatrix}
0 \\
b-Q(\mu,M)\vtight 
\end{bmatrix}.
\end{equation*}
Define $W:=(\Dsqr)\trasp \Dsqr=D_s^{-1}$, and let $\vQ(\cdot)$ be a shorthand for $Q(\cdot,M)\vtight$ (as $M$ is a fixed parameter). Then, for any $w_a, w_b\in \setW$, we can expand \eqref{eq:cocoerc} by using \eqref{eq:TD}, obtaining
\begin{equation}\label{eq:QW}
\begin{split}
(w_a&-w_b)\trasp(\TD(w_a)-\TD(w_b))-\beta\|\TD(w_a)-\TD(w_b)\|^2 \\
={}& (\Dsqr w_a -\Dsqr w_b)\trasp\begin{bmatrix}
F(x_a)-F(x_b) \\
\vQ(\mu_b)-\vQ(\mu_a)
\end{bmatrix}  \\
& \hspace{6em} - \beta\left\|\Dsqr\begin{bmatrix}
F(v_a)-F(v_b) \\
\vQ(\mu_b)-\vQ(\mu_a)
\end{bmatrix}\right\|^2   \\
= &\begin{bmatrix}
x_a-x_b \\
\mu_a-\mu_b
\end{bmatrix}\traspm\begin{bmatrix}
F(x_a)-F(x_b) \\
\vQ(\mu_b)-\vQ(\mu_a) 
\end{bmatrix}  \\
&-\beta\begin{bmatrix}
F(x_a)-F(x_b) \\
\vQ(\mu_b)-\vQ(\mu_a)  
\end{bmatrix}\traspm W \begin{bmatrix}
F(x_a)-F(x_b) \\
\vQ(\mu_b)-\vQ(\mu_a) 
\end{bmatrix},
\end{split}
\end{equation}
for all $y_a,y_b\in\setX\times\setM$, where the last equality follows from the definition of $\setW_j$ and by expanding the norm. Matrix $W$ can be written as $W=\left[\begin{smallmatrix}
W_{11} & W_{12} \\
W_{21} & W_{22}
\end{smallmatrix}\right]$, where $W_{11} \in \mathbb{R}^{nN \times nN}$, $W_{12} \in \mathbb{R}^{nN \times m}$, $W_{33} \in \mathbb{R}^{m \times m}$ are:
\begin{equation*}
\begin{split}
W_{11}&=\tau(I_n-\tau^2A\trasp A)^{-1}, \\
W_{12}&=W_{21}\trasp =\tau^2(I_n-\tau^2A\trasp A)^{-1}A\trasp,\\
W_{22}&=\tau I_m+\tau^3A(I_n-\tau^2A\trasp A)^{-1}A\trasp.
\end{split}
\end{equation*}
Expanding the inner product in \eqref{eq:QW} with respect to the matrix blocks $W_{11}, W_{12} , W_{21} , W_{33}$ we obtain
\begin{equation*}
\begin{split}
	\beta&(F(x_a)-F(x_b))\trasp\big[\tfrac{1}{\beta}(x_a-x_b) \\
	& - W_{11}(F(x_a)-F(x_b)) -2 W_{12}(\vQ(\mu_b)-\vQ(\mu_1))\big]  \\ 
	& +\beta(\vQ(\mu_b)-\vQ(\mu_a))\trasp\big[\tfrac{1}{\beta}(\mu_a-\mu_b) \\
	& \hspace{6em} -W_{22}(\vQ(\mu_b)-\vQ(\mu_a))\big] \\
	={}&(F(x_a)-F(x_b))\trasp(x_a-x_b)  \\ 
	& -\beta(F(x_a)-F(x_b))\trasp W_{11}(F(x_a)-F(x_b)) \\
	&  -2\beta(F(x_a)-F(x_b))\trasp W_{12}(\vQ(\mu_b)-\vQ(\mu_a)) \\ 
	&  +(\vQ(\mu_b)-\vQ(\mu_a))\trasp(\mu_a-\mu_b) \\
	&   -\beta(\vQ(\mu_b)-\vQ(\mu_a))\trasp W_{22}(\vQ(\mu_b)-\vQ(\mu_a)). 
\end{split}
\end{equation*}
Setting $p_\tau:=(I-\tau^2A\trasp A)^{-1/2}(F(x_a)-F(x_b))$ and $q_\tau:=\tau(I-\tau^2A\trasp A)^{-1/2}A\trasp(\vQ(\mu_b)-\vQ(\mu_a))$ above we obtain 
\begin{equation}
\begin{split}
	(&F(x_a)-F(x_b))\trasp(x_a-x_b) \\
	& +(\vQ(\mu_b)-\vQ(\mu_a))\trasp(\mu_a-\mu_b) \\
	& - \beta\tau(\vQ(\mu_b)-\vQ(\mu_a))\trasp(\vQ(\mu_b)-\vQ(\mu_a)) \\
	& -\beta\tau(p_{\tau}+q_{\tau})\trasp(p_{\tau}+q_{\tau}) \\
	&\geq \alpha\|x_a-x_b\|^2 + 2h\zeta\tight \\
	&\hspace{3em} -2\beta\tau h(\tight)^2 - 2\beta\tau(p_{\tau}\trasp p_{\tau}+q_{\tau}\trasp q_{\tau}), \label{eq:H3}
\end{split} 
\end{equation}
where for the last inequality we used, in order, $(i)$ strong monotonicity of $F$ (cf.~Assumption~\ref{assum_monotonicity_continuity}), $(ii)$ Lemma~\ref{lem:Qbound}, $(iii)$ $\|\vQ(\mu_b)-\vQ(\mu_a)\|^2\leq 2h(\tight)^2$---which follows from the same arguments used in the proof of Lemma~\ref{lem:Qbound}---and $(iv)$ $(p_\tau+q_\tau)\trasp(p_\tau+q_\tau) \leq 2(p_{\tau}\trasp p_\tau+q_{\tau}\trasp q_{\tau})$. Expanding the term containing $p_{\tau}, q_{\tau}$ in \eqref{eq:H3} we get
\begin{equation}\label{eq:H4}
	\begin{split}
	&\alpha\|x_a-x_b\|^2 + 2h\zeta\tight -2\beta\tau h(\tight)^2 \\
	&  - 2\beta\tau (F(x_a)\! -\! F(x_b))\trasp (I_n-\tau^2 A\trasp A)^{-1} (F(x_a)\!-\! F(x_b)) \\
	&  - 2\beta\tau^3 (\vQ(\mu_b)-\vQ(\mu_a))\trasp \\
	& \hspace{7em} \cdot A(I_n-\tau^2 A\trasp A)^{-1} A\trasp(\vQ(\mu_b)-\vQ(\mu_a)) \\
	& \overset{(a)}{\geq} \alpha\|x_a-x_b\|^2 + 2h\zeta\tight -2\beta\tau h(\tight)^2 \\
	& \quad - 2\beta\tau \|F(x_a)\! -\! F(x_b)\|^2 \cdot \|(I_n-\tau^2 A\trasp A)^{-1}\| \\
	& \quad - 2\beta\tau^3 \|\vQ(\mu_b)-\vQ(\mu_a)\|^2 \cdot \|(I_n-\tau^2 A\trasp A)^{-1}\|\cdot \|A\|^2 \\
	& \overset{(b)}{\geq} (\alpha-2\beta\tau L_F^2\|(I_n-\tau^2 A\trasp A)^{-1}\|)\|x_a-x_b\|^2 \\
	& \quad + 2h\zeta\tight -2\beta\tau h(\tight)^2 \bigg( 1+ \frac{2\tau^2}{1-\tau^2\|A\|^2} \|A\|^2\bigg),
	\end{split}
\end{equation} 
where $(a)$ is obtained by applying the Cauchy-Schwarz inequality, and in $(b)$ we use the Lipschitz continuity of $F$ (cf.~Assum.~\ref{assum_monotonicity_continuity}), $\|\vQ(\mu_b)-\vQ(\mu_a)\|^2\leq 2h(\tight)^2$, and the triangle inequality. Notice that for the last term in \eqref{eq:H4},
\begin{multline}\label{eq:H4bound}
	2\beta\tau h(\tight)^2 \bigg( 1+ \frac{2\tau^2}{1-\tau^2\|A\|^2} \|A\|^2\bigg) \\ 
	= 	2\beta\tau h(\tight)^2 \frac{1+\tau^2\|A\|^2}{1-\tau^2\|A\|^2}
		\leq 2\beta\tau h(\tight)^2\frac{1+\|A\|^2}{1-\tau^2\|A\|^2}
\end{multline}
holds for any choice of $\tau \in \big(0,  \max\big\{\tfrac{1}{\|A\|},1\big\}\big)$. Recall that by invoking \cite[Thm.~12.5.2]{Pang1}, our objective is to show that \eqref{eq:cocoerc} holds for some $\tau>0$ and $\beta>\tfrac{1}{2}$. Then, by inspecting \eqref{eq:H4} and using \eqref{eq:H4bound}, to achieve this it is sufficient to guarantee
\begin{align*}
	&\alpha-\tau L_F^2\|(I_n-\tau^2 A\trasp A)^{-1}\| > 0,  \\
	&2h\zeta\tight - \tau h(\tight)^2\frac{1+\|A\|^2}{1-\tau^2\|A\|^2} > 0, \; \text{if $1\leq h \leq M$}.
\end{align*}
Solving the quadratic expressions above with respect to $\tau$ results in the admissible range of values in \eqref{eq:tau_cond} (these are also satisfying $\tau \in  \big(0,  \max\big\{\tfrac{1}{\|A\|},1\big\}\big)$, required for \eqref{eq:H4bound} to hold).
Therefore, for any $\tau$ satisfying this condition, $\TD$ is co-coercive with $\beta>\tfrac{1}{2}$ on the entire domain $\setW$, which in turn implies that co-coercivity of $\TD$ holds on each subdomain $\setW_j$, $j=1,\ldots,q$, with the same modulus. By \cite[Thm.~12.5.2]{Pang1}, this is sufficient to guarantee the convergence of \eqref{eq:iter_subset} to a solution of the VI$(\setX\times\setM_j,T)$, thus concluding the proof. \hfill $\blacksquare$ 

\subsection{Proof of Theorem 2}
\begin{figure}
	\centering
	\includegraphics[width=0.8\columnwidth]{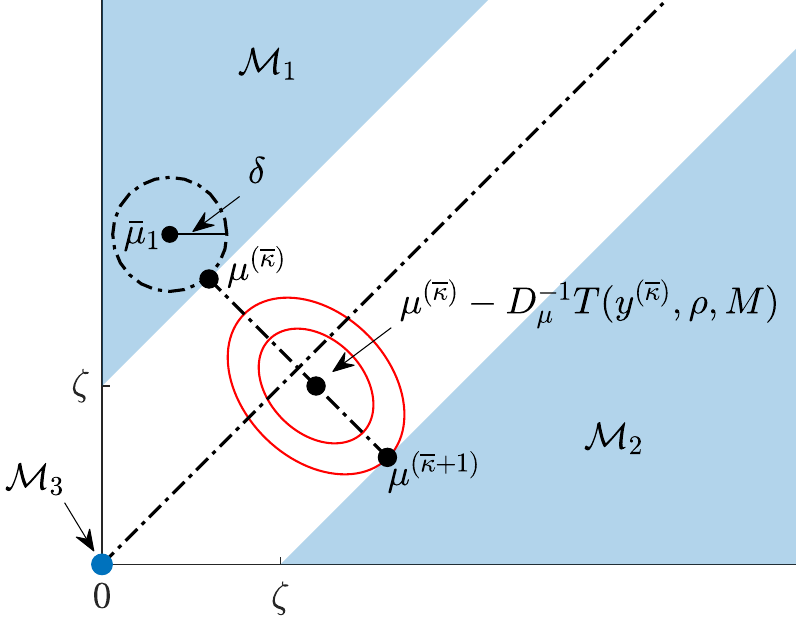}
	\caption{Domain $\setM$ of the Lagrange multipliers associated to the coupling constraints, for the case $m=2$. Notice the minimum distance $\zeta$ between any two subdomains of $\setM$ involves the origin as one of these subdomains.} 
	\label{fig:setM}
\end{figure}

Fix any $\tau$ satisfying the conditions of Lemma~\ref{lem:Tproperties} and \eqref{eq:thm_bound_tau}.
The sequence $\{y^{(\kappa)}\}_{\kappa=1,2,\ldots}$ (where $y^{(\kappa)} = (x^{(\kappa)}, \mu^{(\kappa)})$) generated by Algorithm \ref{Algorithm_apriori} lives in a compact set since $\setX$ and $\setM$ are compact (see Assumption \ref{ass:ass_setM}). As such, by the Bolzano-Weierstrass theorem \cite[Thm. 3.6]{Rudin1976}, there exist convergent subsequences, i.e., the set
\begin{multline}\label{eq:wlimitset}
	\Omega:=\Big\{\bar{y}=(\bar{x},\bar{\mu})\colon\; \exists \text{ subsequence } \{\kappa_i\}_{i\in\mathbb{N}} \\ 
	\text{such that } \lim_{i\to\infty}\kappa_i = \infty,\,\lim_{i\to\infty}y^{(\kappa_i)} = \bar{y} \Big\}, 
\end{multline}
containing the limit points of $\{y^{(\kappa)}\}$ is non-empty; see, e.g., \cite[p. 48]{Rudin1976}. We will show that $\Omega$ is a singleton for any $\tau$ satisfying \eqref{eq:tau_cond}--\eqref{eq:thm_bound_tau}, which implies that the iterates generated by Algorithm \ref{Algorithm_apriori} have a unique limit point, hence they converge.
To achieve this, we assume for the sake of contradiction that there exist two cluster points $\bar{y}_1,\bar{y}_2\in\Omega$, where $\bar{y}_1=(\bar{x}_1,\bar{\mu}_1)$ and $\bar{y}_2=(\bar{x}_2,\bar{\mu}_2)$. Moreover, we assume that $\bar{\mu}_1\in\setM_i$, and $\bar{\mu}_2\in\setM_j$, with $i\neq j$. Note that if this were not the case, then we would be in a trivial case where $\bar{y}_1 = \bar{y}_2$, due to co-coercivity of $T$ (see Lemma~\ref{lem:Tproperties})---by which Algorithm~\ref{Algorithm_apriori} converges to a unique solution when restricted to any convex subdomain $\setX\times\setM_j$, $j=1,\ldots,q$. To ease the notation in the remainder of the proof, we assume without loss of generality that $\bar{\mu}_1\in\setM_1$, $\bar{\mu}_2\in\setM_2$ (see Fig.~\ref{fig:setM}).
By \eqref{eq:wlimitset} there exist an infinite subsequence $\{\kappa_i\}_{i\in\mathbb{N}}$ of the iterates generated by Algorithm \ref{Algorithm_apriori} whose elements get arbitrarily close to $\bar{\mu}_1$ while staying in $\setM_1$ where this cluster point belongs (similarly for $\bar{\mu}_2$). We then have that for any $\delta>0$, there exists $\tilde{\kappa}$ such that for all $\kappa_i\geq \tilde{\kappa}$, $\|y^{(\kappa_i)}-\bar{y}_1\|\leq \delta$; this implies $\|x^{(\kappa_i)}-\bar{x}_1\|\leq \delta$ and $\|\mu^{(\kappa_i)}-\bar{\mu}_1\|\leq \delta$. 

Due to our contradiction hypothesis (recall that $\{\kappa_i\}_{i\in\mathbb{N}}$ is a subsequence), the sequence of iterates generated by Algorithm \ref{Algorithm_apriori} would be leaving $\setM_1$ towards $\setM_2$ infinitely often. Denote then by $\overline{\kappa} > \tilde{\kappa}$ the smallest index of the subsequence such that
$\mu^{(\overline{\kappa})}\in\setM_1$, but $\mu^{(\overline{\kappa}+1)}\in\setM_2$, i.e., after the $\overline{\kappa}$-th iterate the original sequence would jump to $\setM_2$ (for the first time after $\tilde{\kappa}$).
For this jump to occur, the unprojected solution for the Lagrange multipliers must be ``closer'' to $\setM_2$ than to any other sub-domain of $\setM$. 
To see this more formally, let $\Dmu$ denote the lower block-row of $D^{-1} = \left[\begin{smallmatrix} \tau I_{nN} & 0\\ 2A \tau^2& \tau I_m\end{smallmatrix}\right]$, corresponding to the Lagrange multiplier update in line~\ref{alg:update_dyn} of Algorithm \ref{Algorithm_apriori}.
By definition of $\setM$, such a jump requires the Euclidean distance between the unprojected gradient step at $\overline{\kappa}+1$ and  $\mu^{(\overline{\kappa})}$ to satisfy
\begin{equation}
\|\mu^{(\overline{\kappa})}-\Dmu T(y^{(\overline{\kappa})},\rho,M)-\mu^{(\overline{\kappa})}\| > \zeta/2. \label{eq:contra}
\end{equation}
Figure~\ref{fig:setM} illustrates this construction: \eqref{eq:contra} describes the minimum distance for a jump to occur. This is when the ellipsoidal contour levels according to which the projection is performed (skew projection defined by matrix $D$) have their major axis aligned between subdomains as in Figure~\ref{fig:setM} (solid red ellipses). For this two-dimensional example this distance would then be half the width of the white stripe, i.e., $\zeta/\sqrt{2}$. We rather impose $\zeta/2$ (which is smaller) in \eqref{eq:contra}, to account for the case where one of the subdomains is the origin ($\setM_3$). However, 
\begin{align}
\|&\mu^{(\overline{\kappa})}-\Dmu T(y^{(\overline{\kappa})},\rho,M)-\mu^{(\overline{\kappa})}\| \nonumber\\
	&= \tau\|-2\tau A(F(x^{(\overline{\kappa})}) +A\trasp \mu^{(\overline{\kappa})}) \nonumber \\
	&\hspace{0.4cm}+ Ax^{(\overline{\kappa})}-b+Q(\mu^{(\overline{\kappa})},M)\vtight\|  \nonumber\\
        &= \tau\|-2\tau A(F(x^{(\overline{\kappa})}) - F(\bar{x}_1) +A\trasp (\mu^{(\overline{\kappa})} - \bar{\mu}_1)) \nonumber \\
	&\hspace{0.4cm}+ A(x^{(\overline{\kappa})}-\bar{x}_1) + Q(\mu^{(\overline{\kappa})},M)\vtight \nonumber \\
	&\hspace{0.4cm}-2\tau A(F(\bar{x}_1)+A\trasp \bar{\mu}_1) +(A\bar{x}_1 -b) \|  \nonumber\\ 
	&\leq \tau^2\|2A(F(\bar{x}_1)+A\trasp \bar{\mu}_1)\| + \tau\|A\bar{x}_1 -b\| \nonumber \\ 	
	&\hspace{0.4cm} + \tau \|Q(\mu^{(\overline{\kappa})},M)\vtight\| + \tau \|A\| \|x^{(\overline{\kappa})}-\bar{x}_1\| \nonumber\\
	&\hspace{0.4cm}+ 2\tau^2\big(\|A(F(x^{(\overline{\kappa})}) - F(\bar{x}_1))\| +\|AA\trasp(\mu^{(\overline{\kappa})} -  \bar{\mu}_1)\|\big) \nonumber \\
	& \leq (\tau^2+\tau)\bar{R} + \tau\tight\sqrt{m-M} \nonumber \\
	&\hspace{0.4cm}+ \tau\delta \big(2\tau (L_F\|A\| + \|AA\trasp\|) + \|A\|\big), \label{eq:upper_bounding}
	\end{align}
where the first equality follows from the definition of $\Dmu$ and $T$, and the second one by adding and subtracting $F(\bar{x}_1)$, $A\trasp \bar{\mu}_1$ and $A\bar{x}_1$. The first inequality is due to the triangle inequality, while the last one follows from the previous one by upper-bounding (i) the first two terms using the definition of $\bar{R}$; (ii) $\|Q(\mu^{(\overline{\kappa})},M)\vtight\|$ by $\tight\sqrt{m-M}$ based on its definition; and (iii) the last three terms using $\|F(x^{(\overline{\kappa})}) - F(\bar{x}_1)\| \leq L_F \|x^{(\overline{\kappa})}-\bar{x}_1\|$ by Assumption~\ref{assum_monotonicity_continuity}, and $\|x^{(\overline{\kappa})}-\bar{x}_1\|\leq \delta$, $\|\mu^{(\overline{\kappa})}-\bar{\mu}_1\|\leq \delta$. By \eqref{eq:upper_bounding}, and choosing $\tau$ as in \eqref{eq:thm_bound_tau}, we have that 
\begin{align}
\|\mu^{(\overline{\kappa})}-\Dmu T(y^{(\overline{\kappa})},\rho,M)-\mu^{(\overline{\kappa})}\| < \frac{\zeta}{2} + \bar{K}\delta, \label{eq:contra1}
\end{align}
where $\bar{K}$ is a constant, emanating from the coefficient of $\delta$ in \eqref{eq:upper_bounding} when substituting for $\tau$ the upper-bound in \eqref{eq:thm_bound_tau}. 
\KM{Note that $\overline{\kappa}$ is a function of $\delta$, as it depends on $\tilde{\kappa}$, which in turn depends on $\delta$.
Since $\delta$ is arbitrary, taking $\lim \sup_{\delta \to 0}$ in \eqref{eq:contra1} and $\lim \inf_{\delta \to 0}$ in \eqref{eq:contra} leads to 
\begin{align}
\lim \sup_{\delta \to 0}~ &\|\mu^{(\overline{\kappa})}-\Dmu T(y^{(\overline{\kappa})},\rho,M)-\mu^{(\overline{\kappa})}\| < \frac{\zeta}{2},\\
\lim \inf_{\delta \to 0}~ &\|\mu^{(\overline{\kappa})}-\Dmu T(y^{(\overline{\kappa})},\rho,M)-\mu^{(\overline{\kappa})}\| > \frac{\zeta}{2},
\end{align}
establishing a contradiction.}
Then $\bar{\mu}_2,\bar{\mu}_1\in\setM_1$, i.e., all cluster points should be in the same subdomain of $\setM$. As Lemma~\ref{lem:Tproperties} establishes co-coercivity of $T$ on each subdomain $\setX\times\setM_j$, $j=1,\ldots,q$, 
it must be $\bar{\mu}_2=\bar{\mu}_1$, i.e., $\Omega$ is a singleton, implying that Algorithm \ref{Algorithm_apriori} converges. \hfill $\blacksquare$

\subsection{Proof of Theorem \ref{thm:apriori_guarantees}} 

The elements of the minimal compression set $I$ of Algorithm 1 can belong to one or both of the following sets:
\begin{enumerate}
 \item The subset $I_1$ of samples that form a minimal compression for $x^*$. Note that since Algorithm 1 converges to the point $(x^*, \mu^*)$ for a fixed choice of $M$, $Q(\mu^*,M)$ will be a fixed quantity. Then Algorithm 1 will converge to a solution of 
 \begin{align}
 	& \ \text{Find} \ x^* \in \widehat{\Pi}_K \   \text{such that} \nonumber \\ &F( x^*)\trasp(x- x^*) \geq 0 \ \text{for any} \ x \in  \widehat{\Pi}_K, \label{proof_VI}
 \end{align}  
where $\widehat{\Pi}_K$ denotes the polytope obtained from $\Pi_K$ by tightening at most $M$ coupling constraints, as dictated by \eqref{eq:constr_tight} with $Q(\mu^*,M)$. The constraints in \eqref{proof_VI} are equivalent to $F(x^*)\trasp x \geq F(x^*)\trasp x^*$ for all $x \in \widehat{\Pi}_K$. Then, $x^*$ is the minimiser of
\begin{align}
	& \min_{x \in \mathbb{R}^{nN}} F(x^*)\trasp x  \nonumber \\ 
	& \text{subject to }  x \in  \widehat{\Pi}_K, \label{proof_optimization}
\end{align}
which is unique due to Lemma \ref{existence}.
Since the cost function is linear in $x$ and $\widehat{\Pi}_K$ is convex by Assumption \ref{assum_affine}, we obtain a scenario program as in \cite{CampiGaratti2008}. Applying \cite[Prop.~1]{CampiGaratti2008} or \cite[Section III-B]{MargellosEtAl2015} to \eqref{proof_optimization}, we have that $|I_1|\leq nN$, i.e., the cardinality of a minimal compression set for  $x^*$ is bounded by the dimension of the decision vector $nN$.
\item  The subset $I_2$ of samples whose corresponding coupling constraints intersect $\mathbb{B}(x^*, \rho)$. By construction of Algorithm 1 we have that $|I_2|\leq M$.  
\end{enumerate}
As such, we have that $I=I_1 \cup I_2$ is a compression set with cardinality $|I|=|I_1 \cup I_2| \leq |I_1|+|I_2| \leq nN+M$. Then, by Corollary 2 in \cite{MargellosEtAl2015}, \eqref{eq:Theorem3_guarantees} follows.  \hfill $\blacksquare$ 

\bibliographystyle{abbrv}        
\bibliography{ref_consensus,biblio1,biblio3}           

\end{document}